\documentclass[10pt]{article}
\usepackage{amscd,amsmath,amsfonts,amssymb,latexsym}
\usepackage{epsfig}
\usepackage{graphics}
\usepackage[dvips]{color}
\usepackage{pictexwd,dcpic}
\usepackage[vcentermath]{youngtab}
\begin{document}

\newtheorem{prop}{Proposition}
\newtheorem{claim}[prop]{Claim}
\newtheorem{theorem}[prop]{Theorem}
\newtheorem{hypo}[prop]{Hypothesis}
\newtheorem{problem}[prop]{Problem}
\newtheorem{question}[prop]{Question}
\newtheorem{remark}[prop]{Remark}
\newtheorem{lemma}[prop]{Lemma}
\newtheorem{cor}[prop]{Corollary}
\newtheorem{defn}[prop]{Definition}
\newtheorem{ex}[prop]{Example}
\newtheorem{conj}[prop]{Conjecture}
\newtheorem{calc}[prop]{Calculation}
\newcommand{\ca}[1]{{\cal #1}}
\newcommand{\ignore}[1]{}
\newcommand{\f}[2]{{\frac {#1} {#2}}}
\newcommand{\embed}[1]{{#1}^\phi}
\newcommand{\stab}{{\mbox {stab}}}
\newcommand{\perm}{{\mbox {perm}}}
\newcommand{\codim}{{\mbox {codim}}}
\newcommand{\modulo}{{\mbox {mod}\ }}
\newcommand{\E}{\mathbb{E}}
\newcommand{\M}{\mathbb{M}}
\newcommand{\T}{\mathbb{T}}
\newcommand{\C}{\mathbb{C}}
\newcommand{\Z}{\mathbb{Z}}
\newcommand{\Y}{\mathbb{Y}}
\newcommand{\rel}{ \backslash}
\newcommand{\sym}{{\mbox {Sym}}}
\newcommand{\spec}{{\mbox {Spec}}}
\newcommand{\idealof}{{\mbox {ideal}}}
\newcommand{\trace}{{\mbox {trace}}}
\newcommand{\qed}{{\mbox {Q.E.D.}}}
\newcommand{\diff}{{\cal D}}
\newcommand{\invf}{{\cal F}^{-1}}
\newcommand{\invox}{{\cal O}_X^{-1}}
\newcommand{\tableau}{{\mbox {Tab}}}
\newcommand{\zero}{\overline{0}}
\newcommand{\tildae}{\overset{\sim}{e}}
\newcommand{\tilda}[1]{\overset{\sim}{#1}}
\newcommand{\tildaf}{\overset{\sim}{f}}
\newcommand{\arr}{{\cal r}}
\newcommand{\eat}[1]{}
\newcommand{\proof}{\noindent {\em Proof: }}

\title{Quantum deformations of the restriction of 
some $GL_{mn}(\C)$-modules to $GL_m(\C) \times GL_n(\C)$}
\author{
Dedicated to Sri Ramakrishna \\ \\
Bharat Adsul, Milind Sohoni
\footnote{\tt adsul, sohoni@cse.iitb.ac.in}
 \\
Department of Computer Science and Engg., I.I.T. Bombay
 \\   \\ K. V. Subrahmanyam 
\footnote{\tt kv@cmi.ac.in}\\
 Chennai Mathematical Institute}

\maketitle

\begin{abstract} 
In this paper, we consider the restriction of finite dimensional
$GL_{mn} (\C)$-modules to the subgroup ${GL_m (\C )\times GL_n (\C)}$.
In particular, for a Weyl module $V_{\lambda } (\C^{mn})$ of $U_q
(gl_{mn})$ we ask to construct
a representation $W_{\lambda }$ of $U_q (gl_m )\otimes U_q (gl_n )$
such that at $q=1$, the restriction of $V_{\lambda } (\C^{mn})$ to 
$U_1 (gl_m )\otimes U_1 (gl_n )$ matches its action on $W_{\lambda }$ at
$q=1$. 
Thus $W_{\lambda }$ is a $q$-deformation of the module $V_{\lambda }$. 
We achieve this for (i) $\lambda $ consisting of upto two columns and 
$\lambda =(k)$ (i.e., the $Sym^k $ case) for general $m,n$, and 
(ii) all $\lambda $'s for $m=n=2$. This is 
achieved by first constructing a \eat{(commuting)} $U_q (gl_m )\otimes 
U_q (gl_n )$-module $\wedge^k $, a $q$-deformation 
of the simple $GL_{mn} (\C)$-module $\wedge^k (\C^{mn})$. 
We also construct the bi-crystal basis for (i) $\wedge^k $ and show that
it consists of signed subsets, and for (ii) $Sym^k$ and show that it consist
of unordered monomials. 
Next, we develop $U_q (gl_m ) \otimes U_q (gl_n )$-equivariant 
maps $\psi_{a,b} :\wedge^{a+1} 
\otimes \wedge^{b-1} \rightarrow \wedge^a \otimes
\wedge^b$. 
This is used as the building block to construct the general
$W_{\lambda }$ for the cases listed above. 

\end{abstract} 
\section{Introduction}

$GL_N (\C )$ will denote the general linear group of invertible
$N\times N$ complex matrices, and $gl_N(\C)$ its lie algebra.
Consider the group $GL_m (\C) \times GL_n (\C)$ acting on $X$, 
the space of $m\times n$-matrices with complex entries, as follows:
\[ (a,b)\cdot x \rightarrow a \cdot x \cdot b^{T} \]
where $a \in GL_m (\C), \:  b\in GL_n (\C)$ and $x\in X$. Via this
action, we have a homomorphism 
\[\phi :GL_m (\C) \times GL_n (\C) \rightarrow GL_{mn} (\C)\]
For a Weyl module $V_{\lambda }(X)$, via 
$\phi $, we have:
\[ V_{\lambda }(X)=\oplus_{\alpha ,\beta }\: n^{\lambda}_{\alpha ,\beta}
V_{\alpha } (\C^m ) \otimes V_{\beta }(\C^n ) \]
The numbers $n^{\lambda }_{\alpha ,\beta }$ and its properties are
of abiding interest. Even the simplest question of when is
$n^{\lambda}_{\alpha,\beta}>0$ remains unanswered. 

Our own motivation comes from the outstanding problem of P vs. NP, 
and other computational complexity questions in 
theoretical computer science (see \cite{valiant}). More specifically, 
we look at 
the geometric-invariant-theoretic approach to the problem, as
proposed in 
\cite{gct1,gct2}. In this approach, the general subgroup restriction 
problem, i.e., analysing an irreducible representation of a group
$G$ when restricted to a subgroup $H \subseteq G$, is an important
step. An approach to the problem was presented in \cite{gct4}, 
via the dual notion of FRT-algebras (see, e.g., \cite{klimyk}); 
more on this later.

A useful tool in the analysis of representations of the linear groups
$GL_N (\C)$ (henceforth, just $GL_N $), 
has been the {\it quantizations} $U_q (gl_N )$ of the 
enveloping algebra of the lie algebra $gl_N(\C)$, see
\cite{artin,date,drinfeld,jimbo,kashiwara,kashiwara2,lusztigpnas}. The 
representation theory of $U_q (gl_N )$ mimics that of $GL_N$ and has
contributed significantly to the understanding of the diagonal embedding
$GL_N \rightarrow GL_N \times GL_N $, i.e., in the tensor product
of Weyl modules. This is achieved by the Hopf $\Delta : U_q (gl_N ) 
\rightarrow U_q (gl_N ) \otimes U_q (gl_N )$, 
a $q$-deformation of the diagonal
embedding. However, there seems to be no 
quantization of $\phi : GL_m \times GL_n \rightarrow GL_{mn} $, i.e., 
an algebra map (also $\phi$) $U_q (gl_m )\otimes U_q (gl_n ) 
\rightarrow U_q (gl_{mn})$; perhaps none exists \cite{hayashi}. 

On the other hand, we may separately construct embeddings $ U_q (gl_m )
\rightarrow U_q (gl_{mn})$ and $ U_q (gl_n ) \rightarrow
U_q (gl_{mn})$ which correspond to $\phi$ at $q=1$. 
However, the images $ 
(U_q (gl_m ))$ and $(U_q (gl_n ))$ do not commute within 
$U_q (gl_{mn})$. This prevents the standard $U_q(gl_{mn})$-module 
$V_{\lambda } (\C^{mn})$ from becoming a 
$U_q (gl_m ) \otimes U_q (gl_n )$-module.

This paper aims to constructs a $U_q (gl_m )\otimes U_q (gl_n )$-module 
$W_{\lambda }$ with the following properties.
\begin{itemize}
\item $W_{\lambda }$ has a weight structure which matches that of
$V_{\lambda} (\C^{mn})$. Further, there is a weight-preserving 
bijection $W_{\lambda } \rightarrow V_{\lambda} (\C^{mn})$. 
\item The action of $U_q (gl_m )\otimes U_q (gl_n )$ on $W_{\lambda }$ at $q=1$ matches the 
action of $U_1 (gl_m )\otimes U_1 (gl_n ) $ via the embedding 
$\phi :U_1 (gl_m ) \times U_1 (gl_n ) \rightarrow
U_1 (gl_{mn})$
on $V_{\lambda }(\C^{mn})$.
\end{itemize}

We achieve this construction for (i) $\lambda $ with upto two columns, and 
the $\lambda =(k)$ (i.e., the $Sym^k $ case) for general $m$ and $n$, and (ii)
general $\lambda $ for $m=n=2$. We hope to extend these methods for the general
situation.

This construction is done in three steps. We first construct 
$U_q (gl_m )\otimes U_q (gl_n ) $-modules
$W_{\lambda }$ when $V_{\lambda }=\wedge^k (\C^{mn})$, i.e., 
$\lambda $ is a single column shape. Next, we construct 
$U_q (gl_m )\otimes U_q (gl_n ) $-equivariant maps 
\[ \psi_{a,b} : \wedge^{a+1} (\C^{mn}) \otimes \wedge^{b-1}  (\C^{mn}) \rightarrow
\wedge^{a} (\C^{mn}) \otimes 
\wedge^{b} (\C^{mn}) \] 
whose co-kernel is $W_{\lambda }$ when $\lambda $ has two columns.
Finally, the above map gives us straighetning relations which 
yield the construction of general $W_{\lambda }$ when $m=n=2$. Both, the construction
of $\wedge^k (\C^{mn})$ and the map $\psi_{a,b}$ are deformations of 
the usual $U_1 (gl_{mn})$-structures, at $q=1$.

We use the standard model for $U_q (gl_n )$ and its modules consisting of
semi-standard young tableau, see, e.g., \cite{LT}. Thus
a basis for $V_{\lambda } (\C^{mn})$ is identified with $SS(\lambda ,mn)$, i.e., 
semi-standard tableau of shape $\lambda $ with entries in $[mn]$.

In Section \ref{sec:wedge}, we set up notation and then  
construct the 
$U_q (gl_m )\otimes U_q (gl_n )$-modules $\wedge^k $. 
In Section \ref{sec:cb}, we construct a crystal basis for  
$\wedge^k (\C^{mn})$ and show that signed column-tableaus do indeed
form a bi-crystal basis for the $U_q (gl_m )\otimes U_q (gl_n )$-action
thus validating the construction in \cite{danil}. 
Following this, we move towards constructing the abstract module 
$W_{\lambda }$ for 
general $\lambda $. Section \ref{sec:lambda} develops the general line of 
argument and sets up the agenda. Section \ref{sec:psi} proves some elementary
properties of $U_q (gl_m ) \otimes U_q (gl_n )$-modules in the chosen 
basis parametrized by column tableaus. This is used for an explicit 
construction of $\psi_{a,b}$. In Section \ref{sec:sym}, we show that for a 
special choice of 
$\psi_{1,1}, \psi_{2,1}$ and $\psi_{1,2}$, we obtain the $Sym^k $-case. 
We also construct here the crystal base for $Sym^k $. Finally, in 
Section \ref{sec:22}, we consider the $m=n=2$ case and 
show that the above straightening laws yield 
$W_{\lambda }$ for all $\lambda $ (i.e., with upto four rows). 

The construction in this paper has many similarities with that in 
\cite{gct4}. Indeed, our construction of the basic subspaces 
$\wedge^2 (\C^{mn})$ and 
$Sym^2 (\C^{mn})$ of $\C^{mn} \otimes \C^{mn}$ is identical to that 
in \cite{gct4}. There, these subspaces are used to construct the
$R$-matrix and the dual algebra $GL_q (\overline{\C^{mn}})$ and maps
$GL_q (\overline{\C^{mn}})\rightarrow GL_q (\C^m ) \otimes GL_q (\C^n )$. 
The representation theory of $GL_q (\overline{\C^{mn}})$ does not
quite match that of the standard $GL_q (\C^{mn})$ and thus the
construction of $V_{\lambda } (\C^{mn})$ must follow a different route.
Our construction starts with the same $R$-matrix but bypasses the 
construction of $GL_q (\overline{\C^{mn}})$ to arrive directly at
a $GL_q (\C^m ) \otimes GL_q (\C^n )$-structure for $\wedge^k (\C^m )$.
As in \cite{gct4}, we have the ``compactness'' observation, see 
Proposition \ref{prop:transpose}. However, many other structures of 
\cite{gct4} are as yet missing. A point of difference is that even in 
the $m=n=2$ case, we see {\em over-straightening} in \cite{gct4}, while 
here, we do manage to overcome it by a suitable choice of the maps $\psi$.

\section{The $U_q (gl_m) \otimes U_q (gl_n)$ structure for $\wedge^k
(\C^{mn})$} \label{sec:wedge}

To begin, we lift {\em almost verbatim}, the initial parts of 
Section 2 of \cite{LT}.
$U_q (gl_N )$ is the associative algebra over $\C (q)$ generated by
the $4N-2$ symbols $e_i ,f_i ,i=1,\ldots ,N-1$ and $q^{\epsilon_i }, 
q^{-\epsilon_i }$, $i=1,\ldots ,N$ subject to the relations:
\[ q^{\epsilon_i } q^{-\epsilon_i }= q^{-\epsilon_i } q^{\epsilon_i }=
1, \: \: \: [q^{\epsilon_i } ,q^{\epsilon_j }] =0\]
\[ q^{\epsilon_i } e_j q^{-\epsilon_i }= \left\{ \begin{array}{lr}
q e_j & \mbox{for $i=j$} \\
q^{-1} e_j & \mbox{for $i=j+1$} \\
e_j        & \mbox{otherwise} \\ \end{array} \right. \]
\[ q^{\epsilon_i } f_j q^{-\epsilon_i }= \left\{ \begin{array}{lr}
q^{-1} f_j & \mbox{for $i=j$} \\
q f_j & \mbox{for $i=j+1$} \\
f_j        & \mbox{otherwise} \\ \end{array} \right. \]
\[ [e_i ,f_j ]= \delta_{ij} \frac{q^{\epsilon_i } q^{-\epsilon_{i+1}} -
q^{-\epsilon_i }q^{\epsilon_{i+1}}}{q-q^{-1}} \]
\[ [ e_i ,e_j ]=[f_i ,f_j ]=0 \mbox{   for $|i-j|>1$} \]
\[ e_j e_i^2 -(q+q^{-1})e_i e_j e_i +e_i^2 e_j =
f_j f_i^2 -(q+q^{-1})f_i f_j f_i +f_i^2 f_j =0 \mbox{ when $|i-j|=1$} \]

The subalgebra generated by $e_i ,f_i $ and 
\[ q^{h_i }=q^{\epsilon_i }q^{-\epsilon_{i+1}} \: \:\mbox{\hspace*{0.5cm}} 
q^{-h_i }=q^{-\epsilon_i }q^{\epsilon_{i+1}} \mbox{\hspace*{0.5cm}   for $i=1,\ldots ,
N-1$} \]
is denoted by $U_q (sl_N )$. 

The $U_q (gl_N )$ module $V_{(1^k )}$ (henceforth $\wedge^k (\C^N)$) 
is an $\left( \begin{array}{c}
N \\ k \end{array} \right) $-dimensional $C(q)$-vector space with basis
$\{ v_c \}$ indexed by the subsets $c$ of $[N]$ with $k$ elements, 
i.e., by Young Tableau of shape $(1^k )$ with entries in $[N]$. The 
action of $U_q (gl_N )$ on this basis is given by
\[ q^{\epsilon_i } v_c =\left\{ \begin{array}{lr}
v_c & \mbox{  if $i\not \in c$ } \\
qv_c & \mbox{  otherwise} \end{array} \right. \]
\[ e_i v_c =\left\{ \begin{array}{lr}
0 & \mbox{  if $i+1\not \in c$ or $i \in c$ } \\
v_d & \mbox{  otherwise, where $d=c-\{i+1\}+\{i\}$} \end{array} 
\right. \]
\[ f_i v_c =\left\{ \begin{array}{lr}
0 & \mbox{  if $i+1 \in c$ or $i \not \in c$ } \\
v_d & \mbox{  otherwise, where $d=c-\{i\}+\{i+1\}$} \end{array} 
\right. \]

In order to construct more interesting modules, we use the tensor
product operation. Given two $U_q (gl_N )$-modules $M,L$, we can define
a $U_q(gl_N)$-structure on $M \otimes L$ by putting
\[ \begin{array}{rcl}
q^{\epsilon_i } (u \otimes v) &=& q^{\epsilon_i } u \otimes 
q^{\epsilon_i } v \\
e_i (u \otimes v) &=&  e_i u \otimes v + q^{-h_i } u \otimes e_i v \\
f_i (u \otimes v) &=&  f_i u \otimes q^{h_i } v + u \otimes f_i v \\
\end{array} \]
Indeed, the Hopf map $\Delta : U_q (gl_N ) \rightarrow U_q (gl_N ) \otimes 
U_q (gl_N )$:
\[ \Delta q^{\epsilon_i } =q^{\epsilon_i } \otimes q^{\epsilon_i }, 
\Delta e_i =e_i \otimes 1 + q^{-h_i } \otimes e_i , 
\Delta f_i =f_i \otimes q^{h_i} + 1 \otimes f_i \]
is an algebra homomorphism and makes $U_q (gl_N )$ into a bialgebra.

\subsection{Some basic lemmas}

We consider the $U_q (gl_{mn})$-module $\wedge^p (\C^{mn}) $, i.e., 
the homomorphism $U_q (gl_{mn}) \rightarrow 
End_{\C (q)} (\wedge^p (\C^{mn}))$. 
We gather together some lemmas on this particular action. 
\begin{lemma} \label{lemma:V}
On the module $\wedge^p (\C^{mn}) $, we have:
\begin{itemize}
\item $e_i^2 =0$ for all $i$.
\item $e_i e_j e_i =0$ whenever $|j-i|=1$.
\item $e_i f_{i+1}=e_{i+1} f_i =0$ for all $i$.
\end{itemize}
\end{lemma}

We have this important combinatorial lemma:
\begin{lemma} \label{lemma:perm}
Let $\sigma=[\sigma_1 ,\ldots ,\sigma_n ]$ integers such that the 
set $\{ \sigma_1 ,\ldots ,\sigma_n \}=\{ 1,\ldots ,n \}$.
Then, on the module $\wedge^p (\C^{mn})$, for the monomial 
$e_{\sigma }=e_{\sigma_1 } \ldots e_{\sigma_n }$ there exists 
positive integers $k_1 ,\ldots ,k_s $ such that $\sum_i k_i =n$ and 
\[ e_{\sigma} =e_{n-k_s+1} e_{n-k_s +2} \ldots e_n e_{n-k_s -k_{s-1}+1}
e_{n-k_s -k_{s-1}+2} \ldots e_{n-k_s} \ldots e_1 e_2 \ldots e_{k_1} \]
\end{lemma}

\noindent
An important property of the re-ordering is that either (i) the position
of $e_i $ is to the left of position of $e_{i-1}$ or (ii) is {\bf
immediately to the right}.

\begin{ex}
We may verify that:
\[ e_2 e_6 e_7 e_3 e_5 e_1 e_4 =e_6 e_7 e_5 e_2 e_3 e_4 e_1 \]
with $k_1 =1, k_2 =3 ,k_3 =1, k_4 =2$.
\end{ex}

\begin{cor} \label{cor:perm}
Let $\sigma $ be a permutation on the set $\{ i,\ldots ,j\}$
then for the action on $\wedge^p (\C^{mn}) $ we have:
\begin{itemize}
\item if $k<i-1$ or $k>j+1$ then $e_k e_{\sigma }=e_{\sigma } e_k $.
\item if $i\leq k \leq j$ then $e_k e_{\sigma }=e_{\sigma } e_k =0$.
\item if $k<i$ or $k>j$ then $f_k e_{\sigma }=e_{\sigma } f_k $.
\end{itemize}
\end{cor}

For $i<j$, let $E_{i,j}$ denote the term $[e_i , [e_{i+1},[\ldots 
[e_{j-1} ,e_j ]]]$ and $F_{i,j}$ denote $[[[ f_j, f_{j-1} ],\ldots ,f_i
]]$. 

\begin{lemma}
\[
E_{i,j}(v_c) = \left\{ \begin{array}{ll}
(-1)^{|c \cap [i+1,j]|}v_d & \mbox{  if $j+1 \in c$ and $i \not\in c$,
where $d=c-\{j+1\}+\{i\}$} \\
0 & \mbox{otherwise} 
\end{array} 
\right. 
\]
\[
F_{i,j}(v_c) = \left\{ \begin{array}{ll}
(-1)^{|c \cap [i+1,j]|}v_d & \mbox{  if $j+1 \not\in c$ and $i \in c$,
where $d=c-\{i\}+\{j+1\}$} \\
0 & \mbox{otherwise} 
\end{array} 
\right. 
\]

\eat{Further, if $\alpha \neq 0$, then the following holds:
$$\alpha = (-1)^{|c \cap [i+1,j]|}$$
$$ j+1 \in c, i \not\in c $$
$$j+1 \not\in d, i \in d$$
$$c\cap [1,i-1]=d \cap [1,i-1]$$
$$c\cap [i+1,j]=d \cap [i+1,j]$$
$$c\cap [j+2,\ldots]=d\cap [j+2, \ldots]$$}
\end{lemma}
\noindent
{\bf Proof}: 
We provide a detailed proof for $E_{i,j}$. The proof 
for $F_{i,j}$ is similar. 

We prove this by induction on $j-i$.  The base
case is when $j-i=0$. Here, with the convention that $E_{i,i}=e_i$,
the lemma follows from the definition of the operator $e_i$.

For the inductive case (i.e. $i < j$), consider
$E_{i,j}= [e_i, E_{i+1,j}] = e_i E_{i+1,j} - E_{i+1,j} e_i$.
Thus,
$$E_{i,j}(v_c) = e_i E_{i+1,j}(v_c) - E_{i+1,j} e_i(v_c)$$

\eat{By induction,
\[
E_{i+1,j}(v_c) = \left\{ \begin{array}{ll}
(-1)^{|c \cap [i+2,j]|}v_d & \mbox{  if $j+1 \in c$ and $i+1 \not\in c$,
where $d=c-\{j+1\}+\{i+1\}$} \\
0 & \mbox{otherwise} 
\end{array} 
\right. 
\]
}

Suppose that $E_{i+1,j}(v_c) = 0$, so the first-term in
the above expression is zero. Then, by the induction hypothesis,
either $j+1 \not\in c$ or $i+1 \in c$. 

If $i+1 \not\in c$, then $j+1 \not\in c$. Note that in this case, 
$e_i(v_c)=0$. Thus, $E_{i,j}(v_c)=0$ and $j+1 \not\in c$.

If $j+1 \in c$, then $i+1 \in c$. In this case, if $i \in c$, then
$e_i(v_c)=0$ and thus, $E_{i,j}(v_c)=0$ and $i \in c$. Therefore,
we assume that $i \not\in c$ alongwith $j+1 \in c$ and $i+1 \in c$.
So, we have 
\[
e_i(v_c)=v_d \mbox{ where $d=c-\{i+1\}+\{i\}$}
\]
As, $j+1 \in d$ and $i+1 \not\in d$, by induction hypothesis,
\[
E_{i+1,j}(v_d) = 
(-1)^{|d \cap [i+2,j]|}v_e  \mbox{ where $e=d-\{j+1\}+\{i+1\}$} 
\]
Therefore,
\[
\begin{array}{lll}
E_{i,j}(v_c) & = & -E_{i+1,j} e_i(v_c)\\
	   & = & -E_{i+1,j} (v_d)\\
	   & = & -(-1)^{|d \cap [i+2,j]|}v_e\\
	   & = & (-1)^{|c \cap [i+1,j]|}v_e\\
\end{array}
\]
The last equation follows from the fact that $i+1 \in c$ and 
$d=c-\{i+1\}+\{i\}$.
Also, observe that $e=c-\{j+1\}+\{i\}$.

\eat{If $e_i(v_c)=0$ then $E_{i,j}(v_c)=0$
and we are done.
So, we assume that $e_i(v_c)=v_d$. Then
$$ i+1 \in c, i \not\in c $$
$$i+1 \not\in d, i \in d$$
$$c\cap [1,i-1]=d \cap [1,i-1]$$
$$c\cap [i+2,\ldots]=d\cap [i+2, \ldots]$$
Clearly, $E_{i,j}(v_c)= -E_{i+1,j}(v_d)$. By induction,
$E_{i+1,j}(v_d)=\beta v_e$. If
$\beta = 0$, $E_{i,j}(v_c)=0$ and we are done. Otherwise, by induction, 
$$\beta = (-1)^{|d\cap[i+2,j]|}$$
$$ j+1 \in d, i+1 \not\in d $$
$$j+1 \not\in e, i+1 \in e$$
$$d\cap [1,i]=e \cap [1,i]$$
$$d\cap [i+2,j]=e \cap [i+2,j]$$
$$d\cap [j+2,\ldots]=e\cap [j+2, \ldots]$$
Combining the above two observations, we can easily conclude that
$E_{i,j}(v_c) = \alpha v_e$ where
$$\alpha = -\beta = -(-1)^{|d\cap[i+2,j]|} = (-1)^{|c\cap[i+1,j]|}$$
$$ j+1 \in c, i \not\in c $$
$$j+1 \not\in e, i \in e$$
$$c\cap [1,i-1]=e \cap [1,i-1]$$
$$c\cap [i+1,j]=e \cap [i+1,j]$$
$$c\cap [j+2,\ldots]=e\cap [j+2, \ldots]$$}

Now, we consider the case when $E_{i+1,j}(v_c) \neq 0$. Then,
by induction, we have that $j+1 \in c$ and $i+1 \not\in c$. Further, 
\[
E_{i+1,j}(v_c) = 
(-1)^{|c \cap [i+2,j]|}v_d  \mbox{ where $d=c-\{j+1\}+\{i+1\}$} 
\]

Note that, as $i+1 \not \in c$, $e_i(v_c)=0$. Thus, in this
case, 
\[
\begin{array}{lll}
E_{i,j}(v_c) & = & e_i E_{i+1,j}(v_c) - E_{i+1,j} e_i(v_c) \\
	   & = & e_i((-1)^{|c \cap [i+2,j]|}v_d)\\
	   & = & (-1)^{|c \cap [i+2,j]|} e_i(v_d)\\
	   & = & (-1)^{|c \cap [i+1,j]|} e_i(v_d)
\end{array}
\]
The last equality follows from the observation that $i+1 \not\in c$.

If $i \in c$, then $i \in d$ as well and $e_i(v_d)=0$, consequentially
$E_{i,j}(v_c)=0$ as expected.

If $i \not\in c$, then $i \not\in d$ as well. As $i+1 \in d$, we have
$$E_{i,j}(v_c) = (-1)^{|c \cap [i+1,j]|} e_i(v_d)
 = (-1)^{|c \cap [i+1,j]|} v_e$$
where $e=d-\{i+1\}+\{i\}=c-\{j+1\}+\{i\}$.

\eat{$E_{i+1,j}(v_c) = \beta v_d$ where
$$\beta = (-1)^{|c\cap[i+2,j]|}$$
$$ j+1 \in c, i+1 \not\in c $$
$$j+1 \not\in d, i+1 \in d$$
$$c\cap [1,i]=d \cap [1,i]$$
$$c\cap [i+2,j]=d \cap [i+2,j]$$
$$c\cap [j+2,\ldots]=d\cap [j+2, \ldots]$$

Note that, as $i+1 \not \in c$, $e_i(v_c)=0$. Thus, in this case,
$$E_{i,j}(v_c) = e_i E_{i+1,j}(v_c) - E_{i+1,j} e_i(v_c) 
= \beta e_i(v_d)$$ If $e_i(v_d)=0$,
then we are done. So, we assume that $e_i(v_d)=v_e$. Then, we have
the following conditions
$$ i+1 \in d, i \not\in d $$
$$i+1 \not\in e, i \in e$$
$$d\cap [1,i-1]=e \cap [1,i-1]$$
$$d\cap [i+2,\ldots]=e\cap [i+2, \ldots]$$
Combining the above two observations, we can conclude that
$E_{i,j}(v_c) = \alpha v_e$ where
$$\alpha = \beta = (-1)^{|c\cap[i+2,j]|} = (-1)^{|c\cap[i+1,j]|}$$
$$ j+1 \in c, i \not\in c $$
$$j+1 \not\in e, i \in e$$
$$c\cap [1,i-1]=e \cap [1,i-1]$$
$$c\cap [i+1,j]=e \cap [i+1,j]$$
$$c\cap [j+2,\ldots]=e\cap [j+2, \ldots]$$}
\noindent \qed

\begin{lemma} \label{lemma:comm}
For $i,j,i',j'$ , on $\wedge^k (\C^{mn}) $ we have:
\begin{itemize}
\item[(i)] $[E_{i,j} ,E_{i',j'}]=0$ unless either $j'+1=i$ or $j+1=i'$.
\item[(ii)] $[F_{i,j} ,E_{i',j'}]=0$ unless either $j'=j$ or $i'=i$.
\item[(iii)] $E_{i,j} E_{i',j'}=E_{i'j'}E_{ij}=0$ if $i=i'$ or $j=j'$. 
\item[(iv)] $F_{i,j} E_{i',j'}=E_{i'j'}F_{ij}=0$ if $j+1=i'$ or $i=j'+1$. 
\end{itemize}
\end{lemma}

\subsection{Commuting actions on $\wedge^k (\C^{mn})$}

We are now ready to define two actions, that of $U_q (gl_m)$ and 
$U_q (gl_n )$ on $\wedge^p (\C^{mn})$. This will consist of some 
special elements $(E_i^L ,F_i^L ,q^{\epsilon_i^L })$ and $(E_k^R ,F_k^R
,q^{\epsilon_k^R })$ which will implement the action of $U_q (sl_m )$ 
and $U_q (sl_n )$, respectively.

We consider the free $\Z$-module $\E =\oplus_{i=1}^{mn} \Z \epsilon_i $
and define an inner product by extending 
$< \epsilon_i ,\epsilon_j >=\delta_{i,j}$. Define $\kappa_{i,j}
\in \E$ as $\epsilon_i -\epsilon_j $.

We note that:
\begin{lemma}
For $\alpha \in \E$, we have:
\begin{itemize}
\item $e_j q^{\alpha }=q^{<\alpha ,\kappa_{j+1,j}>} q^{\alpha } e_j $.
\item $f_j q^{\alpha }=q^{<\alpha ,\kappa_{j,j+1}>} q^{\alpha } f_j $.
\item $E_{i,j} q^{\alpha }=q^{<\alpha ,\kappa_{j+1,i}>} q^{\alpha } E_{i,j}$.
\end{itemize}
\end{lemma}

Next, we define the {\bf left operators} using:
\[ \begin{array}{rcl} 
B_i^k &=& \sum_{j=0}^{k-2} -h_{jm+i} \\
A_i^k &=& \sum_{j=k}^{n-1} h_{jm+i} \\
\end{array} \]

We define the map $\phi_L : U_q (gl_m ) \rightarrow U_q (gl_{mn})$ 
as:
\[ \begin{array}{rcccl}
q^{\epsilon_i^L} &=& \phi_L (q^{\epsilon_i })&=& \prod_{j=0}^{n-1} 
q^{\epsilon_{jm+i}} \\ 
E^L_i &=& \phi_L (e_i ) &=& e_i +q^{-h_i}e_{m+i} + \ldots
(\prod_{j=0}^{n-2} q^{-h_{jm+i}} ) e_{(n-1)m+i} \\
&& &=& q^{B_i^1} e_i +q^{B_i^2 }e_{m+i} + \ldots
 q^{B_i^n }  e_{(n-1)m+i} \\
F^L_i &=& \phi_L (f_i ) &=& (\prod_{j=1}^{n-1} q^{h_{jm+i}}) f_i + 
\ldots + q^{h_{(n-1)m+i}}f_{(n-2)m+i}+ f_{(n-1)m+i} \\
&&&=& q^{A_i^1}  f_i + 
\ldots + q^{A_i^{n-1}} f_{(n-2)m+i}+ q^{A_i^{n}}f_{(n-1)m+i} \\
\end{array} 
\]

\begin{prop}
The map $\phi_L : U_q (gl_m )\rightarrow U_q (gl_{mn})$ is an algebra
homomorphism.
\end{prop}

\noindent
{\bf Proof}: The embedding of $\phi_L : U_q (gl_m ) \rightarrow U_q (gl_{mn})$ actually 
comes from:
\[ U_q (gl_m ) \stackrel{\Delta }{\longrightarrow} U_q (gl_m )\otimes \ldots U_q (gl_m )
\rightarrow U_q (gl_{mn}) \]
where (i) there are $n$ copies in the tensor-product, and (ii) $\Delta 
$ is the $n$-way Hopf.
This verifies that $\phi_L $ is an algebra map.

We define the {\bf right operators}: 
\begin{defn}
\[ \begin{array}{rcl}
\beta^k_i &=& \sum_{j=i+1}^m \epsilon_{km+j} -\sum_{j=i+1}^m
\epsilon_{(k-1)m+j} \\
\alpha^k_i &=& \sum_{j=1}^{i-1} \epsilon_{i(k-1)m+j} -\sum_{j=1}^{i-1} 
\epsilon_{km+j} 
\end{array}
\]

We define the ``map'' $\phi_R : U_q (gl_n ) \rightarrow U_q (gl_{mn})$ as:
\[ \begin{array}{rcl} 
\phi_R (q^{\epsilon_k^R}) &=& \prod_{i=1}^m q^{\epsilon_{(k-1)m+i}} \\
\phi_R (E_k^R ) &=&  \sum_{i=1}^m q^{\beta^k_i } E_{(k-1)m+i,km+i-1} \\
\phi_R (F_k^R ) &=&  \sum_{i=1}^m q^{\alpha^k_i } F_{(k-1)m+i,km+i-1} \\
\phi_R (h_k^R )&=& \sum_{i=1}^m \epsilon_{(k-1)m+i} -\epsilon_{km+i} 
\end{array} \]
\end{defn}

\noindent{\bf Remark}: $\phi_R $ serves merely to identify a set of 
elements in $U_q (gl_{mn})$ corresponding to the generators of 
$U_q (gl_m)$. 
Thus, while $\phi_L : U_q (gl_m ) \rightarrow U_q (gl_{mn})$ is an 
algebra homomorphism, the corresponding statement for $U_q (gl_n )$ 
is not true. However, as we will show that the composites:
\[ \begin{array}{rcl}
U_q (gl_m ) & \stackrel{\phi_L }{\longrightarrow } & U_q (gl_{mn}) 
\longrightarrow End_{\C (q)} (\wedge^p (\C^{mn})) \\
U_q (gl_n ) & \stackrel{\phi_R }{\longrightarrow } & U_q (gl_{mn}) 
\longrightarrow End_{\C (q)} (\wedge^p (\C^{mn})) \\
\end{array} \]
are {\bf commuting algebra homomorphisms} making $\wedge^p (\C^{mn})$ 
into a $U_q (gl_m ) \otimes U_q (gl_n )$-module. 

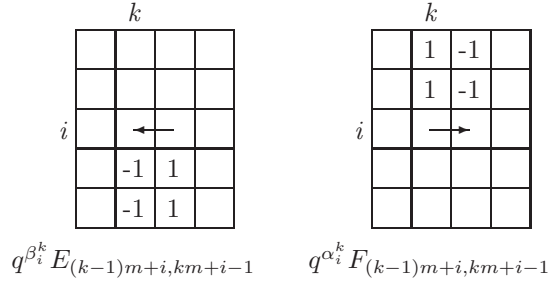
\begin{figure}
\begin{center}
\begin{tabular}{cc}
\begin{picture}(100,100)
\multiput(10,10)(0,15){6}{\line(1,0){60}}
\multiput(10,10)(15,0){5}{\line(0,1){75}}
\put(5,47){\makebox(0,0){\em i}}
\put(32,92){\makebox(0,0){\em k}}
\put(47,47){\vector(-1,0){15}}
\put(32,-2){\makebox(0,0){$q^{\beta_i^k } E_{(k-1)m+i,km+i-1}$}}
\multiput(47,32)(0,-15){2}{\makebox(0,){1}}
\multiput(32,32)(0,-15){2}{\makebox(0,){-1}}
\end{picture}
&
\begin{picture}(100,100)
\multiput(10,10)(0,15){6}{\line(1,0){60}}
\multiput(10,10)(15,0){5}{\line(0,1){75}}
\put(5,47){\makebox(0,0){\em i}}
\put(32,92){\makebox(0,0){\em k}}
\put(32,47){\vector(1,0){15}}
\put(32,-2){\makebox(0,0){$q^{\alpha_i^k } F_{(k-1)m+i,km+i-1}$}}
\multiput(47,62)(0,15){2}{\makebox(0,){-1}}
\multiput(32,62)(0,15){2}{\makebox(0,){1}}
\end{picture}
\end{tabular}
\end{center}
\caption{Terms in the Right Operators} \label{fig:right}
\end{figure}

We will identify $\C^{mn}$ as $\C^m \otimes \C^n $ arranging the typical
element in an $m\times n$ array, reading column-wise from left to right,
and within each column from top to bottom (see below). In this 
notation, see Fig.  \ref{fig:left} for individual terms of the 
left operators and 
Fig. \ref{fig:right} for the right operators. 

\[ \begin{array}{|c|c|c|c|}\hline 
1 & 6 & 11 & 16 \\ \hline 
2 & 7 & 12 & 17 \\ \hline 
3 & 8 & 13 & 18 \\ \hline 
4 & 9 & 14 & 19 \\ \hline 
5 & 10 & 15 & 20 \\ \hline 
\end{array} \]

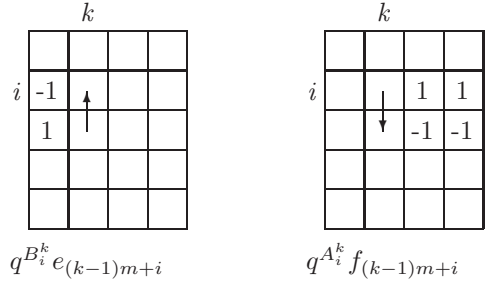
\begin{figure}
\begin{center}
\begin{tabular}{cc}
\begin{picture}(100,100)
\multiput(10,10)(0,15){6}{\line(1,0){60}}
\multiput(10,10)(15,0){5}{\line(0,1){75}}
\put(5,62){\makebox(0,0){\em i}}
\put(32,92){\makebox(0,0){\em k}}
\put(32,47){\vector(0,1){15}}
\put(32,-2){\makebox(0,0){$ q^{B_i^k} e_{(k-1)m+i}$}}
\multiput(17,47)(15,0){1}{\makebox(0,){1}}
\multiput(17,62)(15,0){1}{\makebox(0,){-1}}
\end{picture}
&\begin{picture}(100,100)
\multiput(10,10)(0,15){6}{\line(1,0){60}}
\multiput(10,10)(15,0){5}{\line(0,1){75}}
\put(5,62){\makebox(0,0){\em i}}
\put(32,92){\makebox(0,0){\em k}}
\put(32,62){\vector(0,-1){15}}
\put(32,-2){\makebox(0,0){$ q^{A_i^k } f_{(k-1)m+i}$}}
\multiput(47,47)(15,0){2}{\makebox(0,){-1}}
\multiput(47,62)(15,0){2}{\makebox(0,){1}}
\end{picture}
\end{tabular}
\end{center}
\caption{Terms in the Left Operators} \label{fig:left}
\end{figure}

\subsection{Proofs}

For an operator $O=q^{\mu} E_{i,j}$ (where $\mu \in \E$ is arbitrary) 
let us define $\kappa (O)=\epsilon_{j+1}-\epsilon_i $ and for the 
operator $O=q^{\mu } F_{i,j}$, we define $\kappa (O)$ as 
$\epsilon_i -\epsilon_{j+1}$. We extend this notation so that 
$E_{i,i}=e_i$ (with $\kappa (E_{i,i})=\epsilon_{i+1}-\epsilon_i $) 
and $F_{j,j}=f_j $ (with $\kappa (F_{j,j}) =\epsilon_j -
\epsilon_{j+1}$).

We define ${\cal L}$ and ${\cal R}$ as two sets of operators:
\[ \begin{array}{rcl}
{\cal L} &=&\{ q^{B_i^k } e_{(k-1)m+i} , q^{A_i^k } f_{(k-1)m+i}| 1\leq i \leq
m-1, 1\leq k \leq n \} \\
{\cal R} &=&\{ q^{\beta_i^k } E_{(k-1)m+i, km+i-1} , q^{\alpha_i^k }
F_{(k-1)m+i, km+i-1}| 1\leq i \leq
m, 1\leq k \leq n-1 \} \\
\end{array} \]

Notice that we may write
$E_i^L= \sum_p l_{ip}$ and $E_k^R =\sum_j r_{kj}$ where $l_{ip} \in {\cal
L}$ and $r_{kj} \in {\cal R}$. Whence $[E_i^L ,E_k^R ]$ is expressible
as lie-brackets of elements of ${\cal L}$ and ${\cal R}$. Of course, 
we wish to show that $[E_i^L ,E_k^R ]$ and its three cousins
are actually zero.

\begin{lemma} \label{lemma:reduction}
For any $L\in {\cal L}$ and any $R\in {\cal R}$ if $\langle \kappa
(L),\kappa (R) \rangle \geq 0$ then $[L,R]=0$.
\end{lemma}

\noindent
{\bf Proof}: We first take the case when $\langle \kappa (L),\kappa (R) 
\rangle = 0$. We take for example $L=q^{B_{i'}^{k'}} e_{(k'-1)m+i'}$ 
and $R=q^{\alpha_i^k } F_{(k-1)m+i,km+i-1}$. The condition $\langle 
\kappa (L), \kappa (R) \rangle =0$ implies (see Figs. \ref{fig:left}, 
\ref{fig:right}) that 
\[ \begin{array}{rcl}
F_{(k-1)m+i,km+i-1} q^{B_{i'}^{k'}} &=& q^{B_{i'}^{k'}}
F_{(k-1)m+i,km+i-1} \\
e_{(k'-1)m+i'} q^{\alpha_i^k } &=& q^{\alpha_i^k } e_{(k'-1)m+i'} 
\end{array} \]
Whence 
\[ [L,R]=
q^{B_{i'}^{k'}+\alpha_i^k} [ e_{(k'-1)m+i'} , e_{(k-1)m+i,km+i-1} ]=0 \]
where the last equality follows from Lemma \ref{lemma:comm} (ii).

For the case with $\langle \kappa (L) ,\kappa (R) \rangle =1$, Lemma 
\ref{lemma:comm}, parts (iii),(iv), immediately implies an even stronger claim. \qed

\begin{figure}
\begin{center}
\begin{tabular}{cccc}
\begin{picture}(100,100)
\multiput(10,10)(0,15){3}{\line(1,0){30}}
\multiput(10,10)(15,0){3}{\line(0,1){30}}
\put(5,15){\makebox(0,0)[r]{\em i}}
\put(5,30){\makebox(0,0)[r]{\em i-1}}
\put(32,46){\makebox(0,0){\em k}}
\put(17,46){\makebox(0,0){\em k-1}}
\put(17,17){\vector(1,0){15}}
\put(32,17){\vector(0,1){15}}
\end{picture}
&
\begin{picture}(100,100)
\multiput(10,10)(0,15){3}{\line(1,0){30}}
\multiput(10,10)(15,0){3}{\line(0,1){30}}
\put(5,15){\makebox(0,0)[r]{\em i}}
\put(5,30){\makebox(0,0)[r]{\em i-1}}
\put(32,46){\makebox(0,0){\em k}}
\put(17,46){\makebox(0,0){\em k-1}}
\put(32,17){\vector(0,1){15}}
\put(32,32){\vector(-1,0){15}}
\end{picture}
&
\begin{picture}(100,100)
\multiput(10,10)(0,15){3}{\line(1,0){30}}
\multiput(10,10)(15,0){3}{\line(0,1){30}}
\put(5,15){\makebox(0,0)[r]{\em i}}
\put(5,30){\makebox(0,0)[r]{\em i-1}}
\put(32,46){\makebox(0,0){\em k}}
\put(17,46){\makebox(0,0){\em k-1}}
\put(32,32){\vector(-1,0){15}}
\put(17,32){\vector(0,-1){15}}
\end{picture}
&
\begin{picture}(100,100)
\multiput(10,10)(0,15){3}{\line(1,0){30}}
\multiput(10,10)(15,0){3}{\line(0,1){30}}
\put(5,15){\makebox(0,0)[r]{\em i}}
\put(5,30){\makebox(0,0)[r]{\em i-1}}
\put(32,46){\makebox(0,0){\em k}}
\put(17,46){\makebox(0,0){\em k-1}}
\put(17,32){\vector(0,-1){15}}
\put(17,17){\vector(1,0){15}}
\end{picture}
\\                                               
\begin{picture}(100,100)
\multiput(10,10)(0,15){3}{\line(1,0){30}}
\multiput(10,10)(15,0){3}{\line(0,1){30}}
\put(5,15){\makebox(0,0)[r]{\em i}}
\put(5,30){\makebox(0,0)[r]{\em i-1}}
\put(32,46){\makebox(0,0){\em k}}
\put(17,46){\makebox(0,0){\em k-1}}
\put(32,17){\vector(-1,0){15}}
\put(17,17){\vector(0,1){15}}
\end{picture}
&                                              
\begin{picture}(100,100)
\multiput(10,10)(0,15){3}{\line(1,0){30}}
\multiput(10,10)(15,0){3}{\line(0,1){30}}
\put(5,15){\makebox(0,0)[r]{\em i}}
\put(5,30){\makebox(0,0)[r]{\em i-1}}
\put(32,46){\makebox(0,0){\em k}}
\put(17,46){\makebox(0,0){\em k-1}}
\put(17,17){\vector(0,1){15}}
\put(17,32){\vector(1,0){15}}
\end{picture}
&                                              
\begin{picture}(100,100)
\multiput(10,10)(0,15){3}{\line(1,0){30}}
\multiput(10,10)(15,0){3}{\line(0,1){30}}
\put(5,15){\makebox(0,0)[r]{\em i}}
\put(5,30){\makebox(0,0)[r]{\em i-1}}
\put(32,46){\makebox(0,0){\em k}}
\put(17,46){\makebox(0,0){\em k-1}}
\put(17,32){\vector(1,0){15}}
\put(32,32){\vector(0,-1){15}}
\end{picture}
&                                            
\begin{picture}(100,100)
\multiput(10,10)(0,15){3}{\line(1,0){30}}
\multiput(10,10)(15,0){3}{\line(0,1){30}}
\put(5,15){\makebox(0,0)[r]{\em i}}
\put(5,30){\makebox(0,0)[r]{\em i-1}}
\put(32,46){\makebox(0,0){\em k}}
\put(17,46){\makebox(0,0){\em k-1}}
\put(32,32){\vector(0,-1){15}}
\put(32,17){\vector(-1,0){15}}
\end{picture} \\
\end{tabular}
\end{center}
\caption{The Eight Non-Commuting Terms} \label{fig:noncom}
\end{figure}
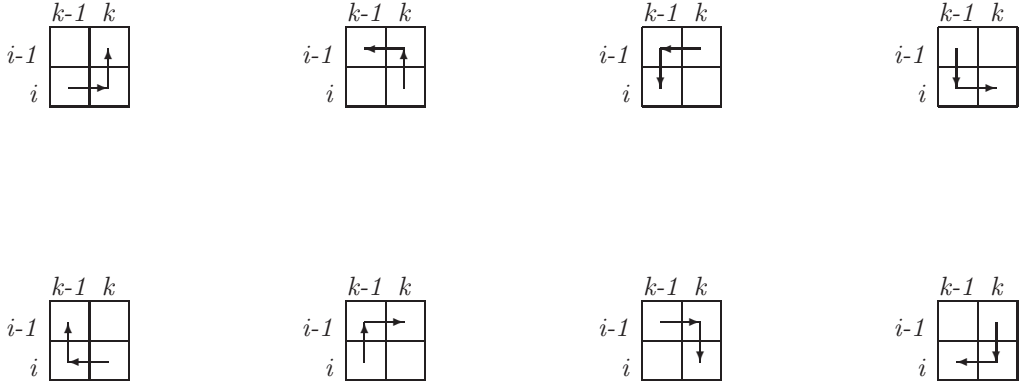

Thus the only non-commuting $(L,R)$ pairs are shown in Fig.
\ref{fig:noncom}.  

By lemma \ref{lemma:reduction}, for the purpose of showing commutation 
we may
as well assume that $n=m=2$. The following argument assumes $n=2$ but
retains $m$ for notational convenience. In other words, we have:
\[ \begin{array}{rcl}
E_i^L &=& e_i + q^{-h_i } e_{m+i} \\
F_i^L &=& q^{h_i} f_i + f_{m+i} \\
\end{array} \]

For $i=1,\ldots ,m$ define $\beta_i , \alpha_i \in \E $ as 
\[ \begin{array}{rcl}
\beta_i &=& \sum_{j=i+1}^m \epsilon_{m+j} -\sum_{j=i+1}^m \epsilon_{j} \\
\alpha_i &=& \sum_{j=1}^{i-1} \epsilon_j -\sum_{j=1}^{i-1} \epsilon_{m+j}
\end{array}
\]
Next, define 
\[ \begin{array}{rcl} 
E^R &=&  \sum_{i=1}^m q^{\beta_i } E_{i,m+i-1} \\
F^R &=&  \sum_{i=1}^m q^{\alpha_i } F_{i,m+i-1} \\
h^R &=& \sum_{i=1}^m \epsilon_i -\epsilon_{m+i} \\
\end{array} \]

Note that $E_1^R =E^R$, $F_1^R =F^R $ and $h_1^R =h^R$.

\begin{lemma}
For $1\leq i \leq m{-1}$,
\begin{itemize}
\item $[e_i ,q^{\beta_{i+1}} E_{i+1,m+i}]=q^{\beta_{i+1}}E_{i,m+i}$. 
\item $[q^{-h_i}e_{m+i} ,q^{\beta_{i}} E_{i,m+i-1}]= 
q q^{\beta_i-h_i}[e_{m+i}, E_{i,m+i-1}]$.
\end{itemize}
\end{lemma}
\proof
We prove the first assertion below. We start with
analyzing
\[
\begin{array}{lll}
[e_i ,q^{\beta_{i+1}} E_{i+1,m+i}] & = &
e_i q^{\beta_{i+1}}E_{i+1,m+i} - q^{\beta_{i+1}}E_{i+1,m+i} e_i \\
& = & 
q^{<\beta_{i+1},-h_i>}
q^{\beta_{i+1}}e_i E_{i+1,m+i} - q^{\beta_{i+1}}E_{i+1,m+i} e_i 
\end{array}
\]

A small calculation shows that $<\beta_{i+1},-h_i>=0$. Therefore,
\[
\begin{array}{lll}
[e_i ,q^{\beta_{i+1}} E_{i+1,m+i}] 
& = & 
q^{\beta_{i+1}}(e_i E_{i+1,m+i} - E_{i+1,m+i} e_i) \\
& = & 
q^{\beta_{i+1}}E_{i,m+i}\\
\end{array}
\]

Now, we turn to the second claim. Towards this, we expand
$[q^{-h_i}e_{m+i} ,q^{\beta_{i}} E_{i,m+i-1}]$ as 
\[
\begin{array}{lll}
& = &
q^{-h_i}e_{m+i}q^{\beta_{i}} E_{i,m+i-1} - q^{\beta_{i}}
E_{i,m+i-1}q^{-h_i}e_{m+i} \\
& = &
q^{-h_i}q^{<\beta_i, -h_{m+i}>}q^{\beta_{i}}e_{m+i} E_{i,m+i-1} - 
q^{\beta_{i}} q^{<-h_i,\kappa_{m+i,i}>}q^{-h_i}E_{i,m+i-1}e_{m+i} \\
\end{array}
\]
We observe that $<\beta_i, -h_{m+i}> = 1$ and 
                $<-h_i,\kappa_{m+i,i}> = 1$. Therefore,
\[
\begin{array}{lll}
[q^{-h_i}e_{m+i} ,q^{\beta_{i}} E_{i,m+i-1}]
& = &
q^{\beta_i-h_i}\left( q e_{m+i} E_{i,m+i-1}
            - q E_{i,m+i-1}e_{m+i} \right) \\
& = &
q. q^{\beta_i-h_i}[e_{m+i}, E_{i,m+i-1}]
\end{array}
\]
\noindent
\qed

\begin{lemma}
$[E^L_i, E^R]=0$
\end{lemma}
\proof
\[
\begin{array}{lll}
[E^L_i, E^R] & = & 
[e_i +q^{-h_i}e_{m+i}, \sum_{j=1}^m q^{\beta_j } E_{j,m+j-1}] \\
& = & 
[e_i ,q^{\beta_{i+1}} E_{i+1,m+i}] + 
[q^{-h_i}e_{m+i} ,q^{\beta_{i}} E_{i,m+i-1}] \\
& = & q^{\beta_{i+1}}E_{i,m+i} + 
q q^{\beta_i-h_i}\left( e_{m+i} E_{i,m+i-1}
            -  E_{i,m+i-1}e_{m+i} \right)
\end{array}
\]
As $\beta_{i} = \beta_{i+1}+\epsilon_{m+i+1}-\epsilon_{i+1}$,
$\beta_{i} - h_i = \beta_{i+1} + \epsilon_{m+i+1} - \epsilon_{i}
= \beta_{i+1} + \kappa_{m+i+1,i}$.

$$[E^L_i, E^R]  = q^{\beta_{i+1}} \left(E_{i,m+i} + 
q.q^{\kappa_{m+i+1,i}}\left( e_{m+i} E_{i,m+i-1}
            -  E_{i,m+i-1}e_{m+i} \right)\right)$$

Now we evaluate the outer bracket at $v_c$. So, we are looking
at $(*)$

$$E_{i,m+i}(v_c) + 
q.q^{\kappa_{m+i+1,i}}\left( e_{m+i} E_{i,m+i-1}(v_c)
            -  E_{i,m+i-1}e_{m+i}(v_c) \right)$$
If $m+i+1 \not\in c$, then all the three terms in the above
expression evaluate to $0$. The middle term certainly evaluates
to $0$ after the application of $e_{m+i}$ (even if
$E_{i,m+i-1}(v_c) \neq 0$).

Similarly, if $i \in c$, then all
the three terms evaluate to $0$. 

So, henceforth, we work with the assumption that
$m+i+1 \in c$ and $i \not\in c$.

Now, we consider the case where $m+i \in c$. In this case,
with $c_1=c-\{m+i+1\}+\{i\}$ and $c_2=c-\{m+i\}+\{i\}$, $(*)$
evaluates to 
\[
\begin{array}{lll}
* & = & (-1)^{|c \cap [i+1,m+i]|}v_{c_1} + 
q.q^{\kappa_{m+i+1,i}}e_{m+i}\left( (-1)^{|c \cap [i+1,m+i-1]|}v_{c_2}
\right) \\
& = & (-1)^{|c \cap [i+1,m+i-1]|}
\left( -v_{c_1} + q.q^{\kappa_{m+i+1,i}} v_{c_1} \right) \\
& = & (-1)^{|c \cap [i+1,m+i-1]|}
\left( -v_{c_1} + q.\frac{1}{q} v_{c_1} \right)\\
& = & 0
\end{array}
\]
Now, we consider the remaining case where $m+i \not\in c$. In this
case, with the notation
$c_1=c-\{m+i+1\}+\{i\}$ and $c_2=c-\{m+i+1\}+\{m+i\}$, $(*)$
evaluates to 
\[
\begin{array}{lll}
* & = & (-1)^{|c \cap [i+1,m+i]|}v_{c_1} -
q.q^{\kappa_{m+i+1,i}}E_{i,m+i-1}(v_{c_2})
 \\
& = & (-1)^{|c \cap [i+1,m+i]|}v_{c_1} -
q.q^{\kappa_{m+i+1,i}}\left((-1)^{|c_2 \cap [i+1,m+i-1]|} v_{c_1}\right)\\

& = & (-1)^{|c \cap [i+1,m+i-1]|}
\left( v_{c_1} - q.q^{\kappa_{m+i+1,i}}v_{c_1}\right)\\
& = & (-1)^{|c \cap [i+1,m+i-1]|}
\left( v_{c_1} - q.\frac{1}{q} v_{c_1}\right)\\
& = & 0
\end{array}
\]

\noindent
\qed

\begin{lemma}
For $1\leq i \leq m-1$,
\begin{itemize}
\item $[f_iq^{h_{m+i}} ,q^{\beta_{i}} E_{i,m+i-1}]=
q q^{h_{m+i}+\beta_i} [f_i, E_{i,m+i-1}]$.
\item $[f_{m+i} ,q^{\beta_{i+1}} E_{i+1,m+i}]=
q^{\beta_{i+1}} [f_{m+i}, E_{i+1,m+i}]$.
\end{itemize}
\end{lemma}
\proof
We start by proving the first claim.
\[
[f_iq^{h_{m+i}} ,q^{\beta_{i}} E_{i,m+i-1}]  = 
f_iq^{h_{m+i}} q^{\beta_{i}} E_{i,m+i-1} -
q^{\beta_{i}} E_{i,m+i-1}f_iq^{h_{m+i}} 
\]
\[
\begin{array}{lll}
f_iq^{h_{m+i}} q^{\beta_{i}} E_{i,m+i-1} & = & 
q^{<h_{m+i}+\beta_i,\kappa_{i,i+1}>}q^{h_{m+i}+\beta_i}
f_iE_{i,m+i-1} \\
& = & q q^{h_{m+i}+\beta_i} f_iE_{i,m+i-1}
\end{array}
\]
\[
\begin{array}{lll}
q^{\beta_{i}} E_{i,m+i-1}f_iq^{h_{m+i}} & = &
q^{\beta_i}E_{i,m+i-1}q^{<h_{m+i}, \kappa_{i,i+1}>}q^{h_{m+i}} f_i \\
& = &
q^{<h_{m+i}, \kappa_{i,i+1}+\kappa_{m+i,i}>}
q^{h_{m+i}+\beta_i}E_{i,m+i-1}f_i \\
& = &
q q^{h_{m+i}+\beta_i}E_{i,m+i-1}f_i 
\end{array}
\]
Thus, 
\[
\begin{array}{lll}
[f_iq^{h_{m+i}} ,q^{\beta_{i}} E_{i,m+i-1}]  & = & 
q q^{h_{m+i}+\beta_i} (f_i E_{i,m+i-1} - E_{i,m+i-1}f_i) \\
& = & q q^{h_{m+i}+\beta_i} [f_i, E_{i,m+i-1}]
\end{array}
\]

Now, we turn to the second claim.

\[
[f_{m+i} ,q^{\beta_{i+1}} E_{i+1,m+i}] = 
f_{m+i} q^{\beta_{i+1}} E_{i+1,m+i} -
q^{\beta_{i+1}} E_{i+1,m+i}f_{m+i} 
\]
\[
\begin{array}{lll}
 f_{m+i} q^{\beta_{i+1}} E_{i+1,m+i} & = & 
q^{<\beta_{i+1},\kappa_{m+i,m+i+1}>}q^{\beta_{i+1}}
f_{m+i}E_{i+1,m+i} \\
& = & q^{\beta_{i+1}} f_{m+i}E_{i+1,m+i}
\end{array}
\]
Thus, 
\[
\begin{array}{lll}
[f_{m+i} ,q^{\beta_{i+1}} E_{i+1,m+i}] & = &
q^{\beta_{i+1}} f_{m+i}E_{i+1,m+i} - 
q^{\beta_{i+1}} E_{i+1,m+i}f_{m+i} 
\\ & = &
q^{\beta_{i+1}} [f_{m+i}, E_{i+1,m+i}]
\end{array}
\]

\begin{lemma}
$[F^L_i, E^R]=0$
\end{lemma}
\proof
\[
\begin{array}{lll}
[F^L_i, E^R] & = & 
[f_iq^{h_{m+i}} +f_{m+i}, \sum_{j=1}^m q^{\beta_j } E_{j,m+j-1}] \\
& = & 
[f_iq^{h_{m+i}},q^{\beta_{i}} E_{i,m+i-1}] +
[f_{m+i} ,q^{\beta_{i+1}} E_{i+1,m+i}]\\ 
& = & 
q q^{h_{m+i}+\beta_i} [f_i, E_{i,m+i-1}] +
q^{\beta_{i+1}} [f_{m+i}, E_{i+1,m+i}]
\end{array}
\]
As $\beta_{i} = \beta_{i+1}+\epsilon_{m+i+1}-\epsilon_{i+1}$,
$\beta_{i} + h_{m+i} = \beta_{i+1} + \epsilon_{m+i} - \epsilon_{i+1}
= \beta_{i+1} + \kappa_{m+i,i+1}$.

$$[E^L_i, E^R]  = q^{\beta_{i+1}} \left(
q q^{\kappa_{m+i,i+1}} [f_i, E_{i,m+i-1}] + [f_{m+i}, E_{i+1,m+i}]
            \right)$$

Now we evaluate the outer bracket at $v_c$. So, we are looking
at $(*)$

$$
q q^{\kappa_{m+i,i+1}} \left(
f_i E_{i,m+i-1}(v_c) - E_{i,m+i-1} f_i(v_c) \right)
+ f_{m+i} E_{i+1,m+i}(v_c) - E_{i+1,m+i}f_{m+i} (v_c)
$$
If $m+i \not\in c$, then all the four terms in the above
expression evaluate to $0$. 
Similarly, if $i+1 \in c$, then all the four terms evaluate to $0$. 

So, henceforth, we work with the assumption that
$m+i \in c$ and $i+1 \not\in c$.

Now, we consider the case where $i \in c$. In this case,
the first term evaluates to 0. If we further assume that
$m+i+1 \not\in c$, then the third term also evaluates to 0.
Overall, 
with $c_1=c-\{i\}+\{i+1\}$ and $c_2=c-\{m+i\}+\{m+i+1\}$, $(*)$
evaluates to 
$$*  =  -q q^{\kappa_{m+i,i+1}} E_{i,m+i-1} (v_{c_1})
- E_{i+1,m+i}(v_{c_2})$$ 
With the notation $d= c-\{m+i\}+\{i+1\}$,
\[
\begin{array}{lll}
* & = & 
-q q^{\kappa_{m+i,i+1}} (-1)^{|c_1 \cap [i+1,m+i-1]|}v_d
- (-1)^{|c_2 \cap [i+2,m+i]|}v_d \\
& = & (-1)^{|c \cap [i+2,m+i-1]|}
\left( qq^{\kappa_{m+i,i+1}} v_d - v_d \right) \\
& = & (-1)^{|c \cap [i+2,m+i-1]|}
\left( q\frac{1}{q} v_d - v_d \right) \\
& = & 0
\end{array}
\]

Now we work with the assumptions $i \in c$ and $m+i+1 \in c$ and
evaluate $(*)$. With these assumptions, the first and the last
term of $(*)$ evaluate to 0. Here, with $c_1=c-\{i\}+\{i+1\}$,
$c_2=c-\{m+i+1\}+\{i+1\}$  and $d=c-\{m+i\}+\{i+1\}$,
$(*)$ evaluates to 
\[
\begin{array}{lll}
* &  = & -q q^{\kappa_{m+i,i+1}} E_{i,m+i-1} (v_{c_1})
+ f_{m+i}((-1)^{|c \cap [i+2,m+i]|} v_{c_2}) \\
& = & 
-q q^{\kappa_{m+i,i+1}} (-1)^{|c_1 \cap [i+1,m+i-1]|}v_d
+ (-1)^{|c \cap [i+2,m+i]|}v_d\\
& = & (-1)^{|c \cap [i+2,m+i-1]|}
\left( qq^{\kappa_{m+i,i+1}} v_d - v_d \right)\\
& = & (-1)^{|c \cap [i+2,m+i-1]|}
\left( q\frac{1}{q} v_d - v_d \right)\\
& = & 0
\end{array}
\]

Now we consider the case with $i \not\in c$. In this case,
the second term in $(*)$ evaluates to 0. As before, if we 
further assume that $m+i+1 \not\in c$, then the third term 
also evaluates to 0.  Overall, 
with $c_1=c-\{m+i\}+\{i\}$, $c_2=c-\{m+i\}+\{m+i+1\}$, 
and $d= c-\{m+i\}+\{i+1\}$,

\[
\begin{array}{lll}
* & = & q q^{\kappa_{m+i,i+1}} f_i((-1)^{|c \cap [i+1,m+i-1]|}v_{c_1})
- E_{i+1,m+i}(v_{c_2}) \\
& = &
q q^{\kappa_{m+i,i+1}} (-1)^{|c \cap [i+1,m+i-1]|}v_d
- (-1)^{|c_2 \cap [i+2,m+i]|}v_d \\
& = & (-1)^{|c \cap [i+2,m+i-1]|}
\left( qq^{\kappa_{m+i,i+1}} v_d - v_d \right) \\
& = & (-1)^{|c \cap [i+2,m+i-1]|}
\left( q\frac{1}{q} v_d - v_d \right) \\
& = & 0
\end{array}
\]

For the only remaining case, we have the assumptions 
$i \not\in c$ and $m+i+1 \in c$.
Here, with $c_1=c-\{m+i\}+\{i\}$,
$c_2=c-\{m+i+1\}+\{i+1\}$  and $d=c-\{m+i\}+\{i+1\}$,
$(*)$ evaluates to 
\[
\begin{array}{lll}
* & = & q q^{\kappa_{m+i,i+1}} f_i((-1)^{|c \cap [i+1,m+i-1]|}v_{c_1})
+ f_{m+i}((-1)^{|c \cap [i+2,m+i]|} v_{c_2}) \\
& = &
q q^{\kappa_{m+i,i+1}} (-1)^{|c \cap [i+1,m+i-1]|}v_d
+ (-1)^{|c \cap [i+2,m+i]|}v_d\\
& = & (-1)^{|c \cap [i+2,m+i-1]|}
\left( qq^{\kappa_{m+i,i+1}} v_d - v_d \right) \\
& = & (-1)^{|c \cap [i+2,m+i-1]|}
\left( q\frac{1}{q} v_d - v_d \right) \\
& = & 0
\end{array}
\]

\noindent
\qed

We have shown that $[E_i^L ,E^R ]=0$ and $[F_i^L ,E^R ]=0$. One can 
similarly show that $[E_i^L ,F^R ]=[F_i^L ,F^R ]=0$. 
We now prepare towards proving $[F_R ,E_R
]=(q^{-h^R}-q^{h^R})/(q-q^{-1})$.

\begin{lemma}
For $i\neq j$, we have:
\[ [q^{\alpha_i } F_{i,m+i-1} ,q^{\beta_j } E_{j,m+j-1}]=0 \]
\end{lemma}

\noindent
{\bf Proof}: 
\[ \begin{array}{rcl}
[q^{\alpha_i } F_{i,m+i-1} ,q^{\beta_j } E_{j,m+j-1}] &=& 
q^{\alpha_i } F_{i,m+i-1} q^{\beta_j } E_{j,m+j-1} -
q^{\beta_j } E_{j,m+j-1} q^{\alpha_i } F_{i,m+i-1} \\
&=& q^{\alpha_i +\beta_j } ( q^{\beta_j (i)-\beta_j (m+i)} F_{i,m+i-1} E_{j,m+j-1} -
 q^{\alpha_i (m+j)-\alpha_i (j}) E_{j,m+j-1} F_{i,m+i-1}) \\
&=& q^a q^{\alpha_i +\beta_j } [F_{i,m+i-1},E_{j,m+j-1}] \\
\end{array} \]
for an appropriate integer $a$ depending on the whether $i\leq j$ or 
not. Now, the only material case for $v_c $ is when $i,m+j \in c$ and 
$j,m+i \not \in c$. We may then verify that $ [F_{i,m+i-1},E_{j,m+j-1}]
v_c =0$. \qed

\begin{lemma}
For $1\leq i \leq m$
\[ [q^{\alpha_i } F_{i,m+i-1} ,q^{\beta_i } E_{i,m+i-1}] =
 q^{\alpha_i +\beta_i } [F_{i,m+i-1},E_{i,m+i-1}] \]
\end{lemma}

\[ \begin{array}{rcl}
[q^{\alpha_i } F_{i,m+i-1} ,q^{\beta_i } E_{i,m+i-1}] &=& 
q^{\alpha_i } F_{i,m+i-1} q^{\beta_i } E_{i,m+i-1} -
q^{\beta_i } E_{i,m+i-1} q^{\alpha_i } F_{i,m+i-1} \\
&=& q^{\alpha_i +\beta_i } ( q^{\beta_i (i)-\beta_i (m+i)} F_{i,m+i-1} E_{i,m+i-1} -
 q^{\alpha_i (m+i)-\alpha_i (i}) E_{i,m+i-1} F_{i,m+i-1}) \\
&=&  q^{\alpha_i +\beta_i } [F_{i,m+i-1},E_{i,m+i-1}] \\
\end{array} \]
This proves the lemma. \qed

Define $\delta_j =\epsilon_j -\epsilon_{m+j}$ and let $v_c \in 
\wedge^p (\C^{mn}) $. 
\begin{lemma}
\[ (q-q^{-1})[F_{i,m+i-1},E_{i,m+i-1}] v_c =
(q^{-\delta_i }-q^{\delta_i} ) v_c \]
\end{lemma}

\noindent
{\bf Proof}: If both $i,m+i \in c$ or both $i,m+i \not \in c$ then 
the equality clearly holds. Now if $i\in c, m+i \not \in c$ then
$q^{\delta_i } v_c = qv_c $ and we have:
\[ \begin{array}{rcl}
(q-q^{-1})[F_{i,m+i-1},E_{i,m+i-1}] v_c &=& (q-q^{-1}) (-v_c )\\
&=& (q^{-\delta_i }-q^{\delta_i} ) v_c \\
\end{array} \]
On the other hand, if $i\not \in c, m+i \in c$, then $q^{\delta_i } v_c =
q^{-1} v_c $ and we have:
\[ \begin{array}{rcl}
(q-q^{-1})[F_{i,m+i-1},E_{i,m+i-1}] v_c &=& (q-q^{-1}) (v_c )\\
&=& (q^{-\delta_i }-q^{\delta_i} ) v_c \\
\end{array} \]
This proves the lemma.

We now prove:
\begin{prop}
Let $h_R =\sum_{i=1}^m \epsilon_i -\epsilon_{m+i}$ then 
\[ [F^R ,E^R ]=\frac{q^{-h_R}-q^{h_R}}{q-q^{-1}} \]
\end{prop}

\noindent
{\bf Proof}: By the above lemmas, we have:

\[ \begin{array}{rcl}
[F^R ,E^R ]&=& \sum_{i=1}^m q^{\alpha_i +\beta_i } [F_{i,m+i-1},E_{i,m+i-1}] \\
\end{array} \]
Whence
\[ \begin{array}{rcl}
(q-q^{-1})[F^R ,E^R ]v_c &=& \sum_{i=1}^m (q^{-\delta_i }-q^{\delta_i
})q^{\alpha_i +\beta_i } v_c  \\
&=& \sum_{i=1}^m q^{\alpha_i +\beta_i -\delta_i } v_c  
-q^{\alpha_i +\beta_i +\delta_i } v_c  \\
\end{array} \]

Now
\[ \alpha_i +\beta_i =(\sum_{j=1}^{i-1} \delta_j )-(\sum_{j=i+1}^m
\delta_j )\]
and thus 
\[ \alpha_i +\beta_i -\delta_i =\alpha_{i-1}+\beta_{i-1}+\delta_{i-1}
= (\sum_{j=1}^{i-1} \delta_j )-(\sum_{j=i}^m
\delta_j )\]
Consequently
\[ \begin{array}{rcl}
(q-q^{-1})[F^R ,E^R ]v_c &=& \sum_{i=1}^m (q^{-\delta_i }-q^{\delta_i
})q^{\alpha_i +\beta_i } v_c  \\
&=& \sum_{i=1}^m q^{\alpha_i +\beta_i -\delta_i } v_c  
-q^{\alpha_i +\beta_i +\delta_i } v_c  \\
&=& (q^{\alpha_1 +\beta_1 -\delta_1 } -q^{\alpha_m +\beta_m
+\delta_m})v_c \\
&=& (q^{-h_R} -q^{h_R }) v_c \\
\end{array} \]
This proves the proposition. \qed

\noindent
We next prove the braid identity. 

\begin{defn}
For $i=1,\ldots ,m$ define $\beta_i , \alpha_i \in \E $ as 
\[ \begin{array}{rcl}
\beta_i &=& \sum_{j=i+1}^m \epsilon_{m+j} -\sum_{j=i+1}^m \epsilon_{j} \\
\beta^*_i &=& \sum_{j=i+1}^m \epsilon_{2m+j} -\sum_{j=i+1}^m \epsilon_{m+j} \\
\end{array}
\]
Next, define 
\[ \begin{array}{rcl} 
E^R &=&  \sum_{i=1}^m q^{\beta_i } E_{i,m+i-1} \\
E^{*R} &=&  \sum_{i=1}^m q^{\beta^*_i } E_{m+i,2m+i-1} \\
\end{array} \]
\end{defn}

Note that $E^{*R}=E_2^R $. We will show that:
\[ (E^R )^2 E^{*R} -(q+q^{-1})E^R E^{*R} E^R + E^{*R} (E^R )^2 =0 \]
We define $g_i =q^{\beta_i } E_{i,m+i-1}$ and $g^*_j =q^{\beta^*_j } 
E_{m+j,2m+j-1}$.

\begin{lemma} \label{lemma:B1}
For distinct $i,j,k\in [m]$ and on $\wedge^p (\C^{mn}) $, we have that 
\[ (g_i g_j +g_j g_i )g^*_k -(q+q^{-1}) (g_i g^*_k g_j +
g_j g^*_k g_i )+g^*_k  (g_i g_j +g_j g_i ) =0 \]
\end{lemma}

\noindent
{\bf Proof}: Let us prove this in several cases. In all cases, we 
will use: 
\[ E_{i,m+i-1} q^{\beta_j } = \left\{ \begin{array}{lr}
q^2 q^{\beta_j } E_{i,m+i-1}  & \mbox{if $i>j$} \\
 q^{\beta_j } E_{i,m+i-1}  & \mbox{if $i\leq j$} \\
 \end{array} \right. \]
\[ E_{i,m+i-1} q^{\beta^*_j } = \left\{ \begin{array}{lr}
q^{-1} q^{\beta^*_j } E_{i,m+i-1}  & \mbox{if $i>j$} \\
 q^{\beta^*_j } E_{i,m+i-1} & \mbox{if $i\leq j$} \\
 \end{array} \right. \]
\[ E_{m+i,2m+i-1} q^{\beta_j } = \left\{ \begin{array}{lr}
q^{-1} q^{\beta_j } E_{i,m+i-1}  & \mbox{if $i>j$} \\
 q^{\beta_j } E_{i,m+i-1} & \mbox{if $i\leq j$} \\
 \end{array} \right. \]

We first consider the case $i<j<k$ and $v_c $ such that $ v
=E_{i,m+i-1}E_{j,j+m-1} E_{m+k,2m+k-1} v_c $, where, by Lemma
\ref{lemma:comm}, the sequence of the operators does not matter. 
Note further that $g_i g_j g^*_k (v_c ) =v^* =q^a \cdot v$. We suppress
the factor $q^a $ uniformly in this proof and in the next lemma 
as well.  We see that:  
\[ \begin{array}{rcl}
(g_i g_j +g_j g_i )g^*_k v_c &=&  
(1 + q^2 ) E_{i,m+i-1}E_{j,j+m-1} E_{m+k,2m+k-1} v_c \\
 &=& (1+q^2 ) v \\
g^*_k (g_i g_j +g_j g_i ) v_c &=&
(q^{-2}+ q^2\cdot q^{-2}) E_{i,m+i-1}E_{j,j+m-1} E_{m+k,2m+k-1} v_c \\
&=& (q^{-2}+1 ) v \\
(g_i g^*_k g_j + g_j g^*_k g_i )v_c &=&
(q^{-1}+  q^{-1}\cdot q^2 ) E_{i,m+i-1}E_{j,j+m-1} E_{m+k,2m+k-1} v_c \\
&=& (q^{-1}+q) v \\
\end{array}
\]
This proves the assertion for $i<j<k$.

Next, let us consider $i<k<j$:
\[ \begin{array}{rcl}
(g_i g_j +g_j g_i )g^*_k v_c &=& (q^{-1}+q ) v \\
g^*_k (g_i g_j +g_j g_i ) v_c &=& (q^{-1}+q ) v \\
(g_i g^*_k g_j + g_j g^*_k g_i )v_c &=& 2 v \\
\end{array}
\]
This proves the assertion for $i<k<j$.

Next, let us consider $k<i<j$:
\[ \begin{array}{rcl}
(g_i g_j +g_j g_i )g^*_k v_c &=& (1+q^{-2} ) v \\
g^*_k (g_i g_j +g_j g_i ) v_c &=& (1+q^{2} ) v \\
(g_i g^*_k g_j + g_j g^*_k g_i )v_c &=& (q+q^{-1}) v \\
\end{array}
\]
This proves the assertion for $k<i<j$ and completes the proof of the 
lemma. \qed

\begin{lemma} \label{lemma:B2}
For distinct $i,j\in [m]$ and on $\wedge^p (\C^{mn}) $, we have that 
\[ (g_i g_j +g_j g_i )g^*_i -(q+q^{-1}) (g_i g^*_i g_j +
g_j g^*_i g_i )+g^*_i  (g_i g_j +g_j g_i ) =0 \]
\end{lemma}

\noindent
{\bf Proof}: There are two cases to consider, viz., $g_i g_i^* v_c =0$ 
and $g_i^* g_i v_c =0$. Let us consider the first case, i.e., 
$g_i g^*_i v_c =0$, in which case we need to show: 
\[  -(q+q^{-1})  g_j g^*_i g_i +g^*_i  (g_i g_j +g_j g_i ) =0 \]
Let $v $ be such that $ E_{m+i,2m+i-1} E_{i,m+i-1}E_{j,m+j-1}v_c
=v $ (see comment in proof of Lemma \ref{lemma:B1}). 
We see that for $j>i$:
\[ \begin{array}{rcl}
g^*_i (g_i g_j +g_j g_i ) v_c &=&  (1+q^2 )v \\
g_j g^*_i g_i v_c &=& q v \\
\end{array}
\]
This proves the lemma for $j>i$. Next, for $j<i$, with
$v = E_{j,m+j-1}E_{m+i,2m+i-1} E_{i,m+i-1}v_c$ and we have:
\[ \begin{array}{rcl}
g^*_i (g_i g_j +g_j g_i ) v_c &=&  (q+q^{-1} )v \\
g_j g^*_i g_i v_c &=& v \\
\end{array} \]

This proves the case when $g_i g^*_i v_c =0$. The other case is 
similarly proved. \qed

\begin{prop}
For $E^R =E_1^R $ and $E^{*R}=E_2^R$, we have:
\[ (E^R )^2 E^{*R} -(q+q^{-1})E^R E^{*R} E^R + E^{*R} (E^R )^2 =0 \]
\end{prop}

\noindent
{\bf Proof}: 
Let 
\[ B=(E^R )^2 E^{*R} -(q+q^{-1})E^R E^{*R} E^R + E^{*R} (E^R )^2 \]
For a given $v_c $, we look at $B\cdot v_c $ and classify the result
by the $U_q (gl_{mn})$ weight. We see that the allowed weights are
$wt(v_c )-\kappa_{m+i,i}-\kappa_{m+j,j}-\kappa_{m+k,k}$ for various 
$i,j,k$. Further, we see that:
\[ \begin{array}{rcl}
   E^R &=& \sum_{i=1}^m g_i \\
   E^{*R} &=& \sum_{i=1}^m g^*_i \\ \end{array} \]
is a separation of $E^R $ and $E^{*R}$ by $U_q (gl_{mn})$-weights. 
Therefore showing $B\cdot v_c =0$ amounts to various cases on $i,j,k$. 
The main cases are settled by Lemmas \ref{lemma:B1},
\ref{lemma:B2}. Other cases are easier. \qed

\begin{prop} 
The map $\phi_R : U_q (gl_n ) \rightarrow End_{\C (q) } (\wedge^p
(\C^{mn})) $ is an algebra homomorphism. At $q=1$, $\phi_R $ factorizes
through $U_q (gl_{mn})$, i.e., 
\[ \phi_R (1): U_1 (gl_n ) \rightarrow U_1 (gl_{mn}) \rightarrow
End_{\C} (\wedge^p (\C^{mn})) \]
\end{prop}

\noindent
The proof is obvious. The family $\{ E_k^R ,F_k^R ,q^{\epsilon_k^R}  \}$ satisfy
all the properties for $U_q (gl_n )$. 
Also note that at $q=1$, $\phi_R (1)$ reduces to the standard injection 
which commutes with $\phi_L (1)$.

\section{The crystal basis for $\wedge^K $} \label{sec:cb}

In this section we examine the crystal structure (see 
\cite{kashiwara,kashiwara2}) of the $U_q (gl_m )\otimes 
U_q (gl_n )$-module $\wedge^K (\C^{m\times n})$. We show that there is 
a sign function $sign^* $ on $K$-subsets of $[mn]$ such that the 
collection ${\cal B}^* =\{ sign^* (c) \cdot v_c \}_c $ is a crystal
basis for $\wedge^K $.

We identify $[mn]$ with
$[m] \times [n]$ and also order the elements as follows:
\[ (1,1) \prec (2,1) \prec \ldots (m,1) \prec (1,2) \prec \ldots 
(m-1,n) \prec (m,n) \] 
In other words $(i,j) \prec (i',j')$ iff either $j<j'$ or 
$j=j'$ with $i<i'$. For $(i,j) \prec (i',j')$, we denote by 
$[ (i,j),(i',j')]$ as the indices between $(i,j)$ and $(i',j')$
including both $(i,j)$ and $(i',j')$.

Recall that (cf. Section \ref{sec:wedge}), 
as a $\C (q)$-vector space, $\wedge^K(C^{mn})$ is 
generated by the basis 
vectors ${\cal B}= \{ v_c | c\subseteq [mn], |c|=K \}$. 
Let us fix an index $i$ and look at the 
sub-algebra $U_i^L $ of $U_q (gl_m )$ generated by 
$E_i^L ,F_i^L $ and $h_i^L $. We define the standard $U_q (sl_2 )$
generated by symbols $e,f,h$ satisfying the following equations:
\[ \begin{array}{rcl}
q^{h}q^{-h}&=& 1 \\
q^h eq^{-h} &=& q^2 e \\
q^h fq^{-h} &=& q^{-2} f \\
ef-fe&=& \frac{e^h -e^{-h}}{q-q^{-1} }\\
\end{array} \]
We use the Hopf $\Delta $:
\[ \Delta q^{h} =q^{h} \otimes q^{h},
  \Delta e =e \otimes 1 + q^{-h } \otimes e ,
   \Delta f =f \otimes q^{h} + 1 \otimes f \]

In other words, they satisfy exactly the same relations that 
$e_i^L ,f_i^L ,h_i^L $ satisfy, including the Hopf.
Clearly, $U_i^L $ is isomorphic to $U_q (sl_2 ) $ 
as algebras
and we denote this isomorphism by $L: U_i^L  \rightarrow U_q (sl_2 )$. 

We construct the $U_q (sl_2 )$-module $\C^2 $ with basis $x_1 ,x_2 $
with the action:
\[ ex_2 =x_1 , ex_1 =0, fx_2 =0, fx_1 =x_2 , q^h x_1 =qx_1 , q^h x_2
=q^{-1} x_2 \]  

With the Hopf $\Delta $ above, $M=\otimes_{i=1}^N \C^2 $ 
is a $U_q (sl_2 )$-module with the basis ${\cal S}=\{ y_1 \otimes 
\ldots \otimes 
y_N | y_i \in \{ x_1 ,x_ 2 \} \} $, and with the action:
\[ e(y_1 \otimes \ldots y_N )=
\sum_j (\prod_{k=1}^{j-1} q^{-h} (y_k ) ) \cdot y_1 \otimes \ldots 
\otimes y_{j-1} \otimes 
e(y_j ) \otimes y_{j+1} \otimes\ldots \otimes y_N \]
A similar expression may be written for the action of $f$. 

Let us identify $[mn]$ with $[m]\times [n]$
and define the {\bf signature} $\sigma_i^L (c)$, for $c \subseteq [mn]$. 
Towards this, we define 
\[ \begin{array}{rcl}
I(c)&=& \{ 1\leq j \leq n \mid \mbox{~both~} (i,j),(i+1,j) \in c \} \\
J(c)&=& \{ 1\leq j \leq n \mid \mbox{~both~} (i,j),(i+1,j) \not \in c \} \\
S(c)&=& \{ (i',j') \in c \mid i'\neq i \mbox{~and~} i' \neq i+1 \}
\end{array} \]
The signature $\sigma_i^L  (c)$ is the tuple $(I(c),J(c),S(c))$. 

Next, for a $\sigma =(I,J,S)$, we define the vector 
space $V_{\sigma ,i}^L $ as the $\C(q)$-span of all elements 
\[ {\cal B}_{\sigma,i}^L=\{ v_c \mid \sigma_i^L (c)=\sigma \} \] 

Let $N=n-|I|-|J|$ and let $M=\otimes^N \C^2 $ be the 
$U_q (sl_2 )$-module as above. 

We prove the following:
\begin{prop} \label{prop:Liso}
Given $\sigma =(I_{\sigma } ,J_{\sigma }, S_{\sigma })$ as above,
\begin{itemize}
\item[(i)] $V_{\sigma ,i }^L$ is a $U_i^L $-invariant subspace.
\item[(ii)] The $U_q (sl_2 )$ module $M$ is isomorphic to the $U_i^L $-module 
$V_{\sigma ,i}^L $ via the isomorphism $L$ above.
\end{itemize}
\end{prop}

\noindent
{\bf Proof}: For any $v_c \in {\cal B}_{\sigma ,i}^L$, if $E_i^L (v_c )=
\sum \alpha (c') \cdot v_{c'}$, then it is clear that $v_{c'} 
\in {\cal B}_{\sigma ,i}^L$ as well. The same holds for $F_i^L $ and 
$h_i^L $. This proves (i) above. 
For (ii), first note that 
\[ E_i^L = \sum_j (\prod_{k=1}^{j-1} q^{-h_{(k-1)m+i}} ) e_{(j-1)m+i} \]
which matches the Hopf $\Delta $ of $U_q (sl_2 )$. 
Next, if $j\in I(c)\cup J(c)$ then the index $j$ is irrelevant to the 
action of $E_i^L $ on $v_c $, whence in the restriction to 
$V_{\sigma ,i}^L $, the indices in $I_{\sigma } \cup J_{\sigma } $ do 
not play a role.

Next, note that $|{\cal B}_{\sigma ,i}^L|=2^N$. 
Assume for simplicity that $I_\sigma \cup J_\sigma =
\{ N+1, \ldots ,n \}$.
Indeed, we may set up a $U_q (sl_2 )$-module isomorphism $\iota_L $ by 
setting 
\[ \iota_L (v_c )=y_1 \otimes \ldots \otimes y_N  \mbox{ such that }
y_k =\left\{ \begin{array}{l}
x_1  \mbox{ iff } (i,k)\in c \\
x_2 \mbox{ otherwise } \end{array} \right. \] 
One may verify that $\iota_L :V_{\sigma ,i}^L \rightarrow M$ 
is indeed equivariant via $L$.  \qed

\begin{prop} \label{prop:Lcb}
The elements ${\cal B}$ is a crystal basis for $\wedge^K (\C^{mn})$ for the 
action of $U_q (gl_m )$. 
\end{prop}

\noindent
{\bf Proof}: This is obtained by first noting that ${\cal S}$ is indeed
a crystal basis for $M$, see \cite{kashiwara}, for example. Next, the 
equivariance of $\iota_L $ shows that for $v_c \in 
{\cal B}_{\sigma ,i}^L $, 
\[ \stackrel{\sim}{E_{i}^L } (v_c )=  \iota_L^{-1} (\: \stackrel{\sim }
{e} (\iota_L (v_c ))) \]
This proves
that ${\cal B}_{\sigma ,i}^L$ is indeed a crystal basis for $V_{\sigma 
,i}^L$. Next, by applying Proposition \ref{prop:Liso} for all $i$ and 
all $\sigma $, we see that $\{ {\cal B}_{\sigma ,i}^{\sigma } | i,
\sigma \} $ together cover ${\cal B}$. \qed

We now move to the trickier $U_q (gl_n)$-action. 
Let us denote by $\epsilon_{i,j}$ the weight $\epsilon_{(j-1)*m+i}$
and $h_{i,j}=\epsilon_{i,j}-\epsilon_{i,j+1}$.
There are two sources of complications. 
\begin{itemize}
\item The operator $E_k^R $ may be re-written as: 
\[  
E_k^R = \sum_i \prod_{a=i+1}^m (q^{-h_{a,k}}) E_{(k-1)m+i,km+i-1}
= \sum_i E_{(k-1)m+i,km+i-1} (\prod_{a=i+1}^m q^{-h_{a,k}}) 
\]
Thus, the Hopf works from the ``right''.
\item For a general $v_c $, if $E_{(k-1)m+i, km+i-1} v_c 
$ is non-zero then it is $\pm v_d $, where $v_d =v_c -(i,k+1)+(i,k)$
where the sign is $(-1)^M$ where $M$ is the number 
of elements in $c\cap [(i+1,k),\ldots ,(i-1,k+1)]$. 
\end{itemize} 

To fix the sign, we first define an ``intermediate global'' sign as 
follows. For 
a set $c \subset [m]\times [n]$, we define $c^* \subset [m]\times 
[n]$ as that obtained by moving the elements
of $c$ to the right, as far as they can go (see 
Example \ref{ex:wedge2}). Note that $F_k^R (c^* )=0$ 
for all $k$ and thus $c^* $ is {\em one of the} lowest weight vectors
in $\wedge^K (\C^{mn})$. 
For an $(i,j) \in c $, 
let $(i ,j^* )$ be its final position in $c^* $. We may define $j^* $
explicitly as $n-| \{ j' | (i,j') \in c , j'>j \}|$. Next, we define for 
$(i,j) \in c$, 
\[ \begin{array}{rcl}
S_{i,j} (c)&=& \{ (i',j') \in c \mid  
(i,j) \prec (i',j') \prec (i',j'^*) \prec (i,j^* ) \} \\
n_{ij} &=& |S_{i,j} (c)| 
\end{array}
\]
Setting $N_c =\sum_{(i,j) \in c } n_{ij}$ we finally define:
\[ \begin{array}{rcl}
sign (c)&=& (-1)^{N_c } \\
sign (d/c )&=& sign(d)/sign(c)  \end{array} \]

\begin{lemma} \label{lemma:ell}
Let $v_c \in {\cal B}_{\sigma ,k}^R $ be such that 
$E_{(k-1)m+i, km+i-1} v_c \neq 0$, then 
\[ E_{(k-1)m+i, km+i-1} v_c  =sign(d/c) v_d \] 
where $v_d =v_c -(i,k+1)+(i,k)$.
\end{lemma}

\noindent
{\bf Proof}: It is clear that $c^* =d^* $ and thus for $(i,k+1) \in c$ 
and $(i,k) \in d$, let $(i,k^* )$ be the final position of both 
$(i,k+1) \in c$ and $(i,k) \in d$. For $(i,k+1) \prec (i',j')$ or
$(i',j') \prec (i,k)$ we have (i) $S_{i',j'} (c)=S_{i',j'} (d)$ and (ii)
$(i',j') \in S_{i,k+1}(c) $ iff $(i',j') \in S_{i,k} (d)$. 

Next, it is clear that (i) $S_{i,k} (d) \supseteq S_{i,k+1} (c)$, and 
(ii) for $(i,k) \prec (i',j') \prec (i,k+1)$, $S_{i',j'} (d) 
\subseteq S_{i',j'} (c)$ and in fact, $S_{i',j'} (c) -S_{i',j'} (d)$ 
can atmost be the element $(i,k+1)$.

Now let us look at $S_{i,k} (d)-S_{i,k+1} (c)$. These contain 
all $(i',j')\in c$ such that 
\[ (i,k )\prec (i',j') \prec (i,k+1) \prec 
(i',j'^* ) \prec (i,k^* ) \]
On the other hand, for $(i',j') \in c$ such that 
$(i,k) \prec (i' ,j') \prec (i,k+1)$, which 
are not counted above, it must be that $(i,k^* ) \prec (i',j'^* ) $
in which case, $S_{i',j'} (c)=S_{i',j'} (d)\cup \{ (i,k+1) \}$. 

In short, for
every $(i',j') \in c$ such that $(i,k) \prec (i',j') \prec (i,k+1)$ 
either it contributes to an increment in $S_{i,k} (d)$ over $S_{i,k+1}
(c)$ or a decrement
in $S_{i',j'} (d)$ over $S_{i',j'} (c)$. Ofcourse, the two cases
are exclusive.

Thus we have 
$sign(d)/sign(c)=(-1)^M $ where $M$ is exactly the number 
of elements in $c\cap [(i+1,k),\ldots ,(i-1,k+1)]$. \qed

Next, we define  a new Hopf $\Delta '$ on $U_q (sl_2 )$ as
\[ \Delta 'q^{h} =q^{h} \otimes q^{h},
  \Delta 'e =1 \otimes e + e \otimes q^{-h} ,
   \Delta 'f =q^h  \otimes f + f \otimes 1 \]
We denote by $M'$, the $U_q (sl_2 )$-module $\otimes^N \C^2 $ via
the Hopf $\Delta '$ and with the basis ${\cal S}=\{ y_1 \otimes 
\ldots \otimes y_N | y_i \in \{ x_1 ,x_ 2 \} \} $. Under $\Delta '$
we have: 
\[ e(y_1 \otimes \ldots y_N )=
\sum_j (\prod_{k=j+1}^{N} q^{-h} (y_k ) ) \cdot y_1 \otimes \ldots 
\otimes y_{j-1} \otimes 
e(y_j ) \otimes y_{j+1} \otimes\ldots \otimes y_N \]

We denote by $U_k^R $ the algebra generated by $E_k^R ,F_k^R , 
h_k^R $ and let $R: U_k^R \rightarrow U_q (sl_2 )$ be the natural
isomorphism. 

As before, we define $\sigma_k^R (c)$ analogously as
\[ \begin{array}{rcl}
I(c)&=& \{ 1\leq i \leq m \mid \mbox{~both~} (i,k),(i,k+1) \in c \} \\
J(c)&=& \{ 1\leq i \leq m \mid \mbox{~both~} (i,k),(i,k+1) \not \in c \} \\
S(c)&=& \{ (i',k') \in c \mid k'\neq k \mbox{~and~} k' \neq k+1 \}
\end{array} \]
Next, for a $\sigma =(I,J,S)$, we define the vector 
space $V_{\sigma ,k}^R $ as the $\C(q)$-span of all elements 
\[ {\cal B}_{\sigma,k}^R=\{ v_c \mid \sigma_k^R (c)=
\sigma \} \] 

Again, as before, let $N=n-|I|-|J|$. \eat{
and let $M=\otimes^N \C^2 $ be the 
$U_q (sl_2 )$-module as above.} Let us also assume, for simplicity
that $I\cup J=\{ N+1, \ldots ,m \} $.

\begin{prop} \label{prop:Riso}
Given $\sigma $ as above,
\begin{itemize}
\item[(i)] $V_{\sigma ,k }^R$ is a $U_k^R $-invariant subspace.
\item[(ii)]The $U_q (sl_2 )$ module $M'$ is isomorphic to the $U_k^R $-module 
$V_{\sigma ,k}^R $ via the isomorphism $R$ above.
\end{itemize}
\end{prop}

\noindent
{\bf Proof}: Part (i) above is obvious. 
For (ii), note that 
\[  
E_k^R = \sum_i E_{(k-1)m+i,km+i-1} (\prod_{a=i+1}^m q^{-h_{a,k}}) \]
which matches the Hopf $\Delta' $ of $U_q (sl_2 )$. 
Again, if $j\in I(c)\cup J(c)$ then the index $j$ is irrelevant to the 
action of $E_k^R $ on $v_c $, whence in the restriction to 
$V_{\sigma ,k}^R $, the indices in $I \cup J $ do 
not play a role.

Next, note that $|{\cal B}_{\sigma ,k}^R|=2^N$. 
Recall that, we have assumed  
that $I\cup J =\{ N+1, \ldots ,m \}$.
Indeed, we may set up a $U_q (sl_2 )$-module isomorphism $\iota_R $ by 
setting 
\[ \iota_R (v_c )= sign(c) \cdot y_1 \otimes \ldots \otimes y_N  \mbox{ such that }
y_i =\left\{ \begin{array}{l}
x_1  \mbox{ iff } (i,k)\in c \\
x_2 \mbox{ otherwise } \end{array} \right. \] 
One may verify (using Lemma~\ref{lemma:ell}) that $\iota_R :V_{\sigma ,k}^R \rightarrow M'$ 
is indeed equivariant via $R$.  \qed

\begin{prop} \label{prop:Rcb}
Let ${\cal B}'=\{ sign(b) \cdot v_b |b \in {\cal B} \}$ be ``signed''
elements. Then the elements ${\cal B}'$ is a crystal basis for 
$\wedge^K (\C^{mn})$ for the action of $U_q (gl_n )$. In other 
words $\stackrel{\sim }{E_k^R} (v_c ) =\pm v_d \cup 0 $.
\end{prop}

\noindent
{\bf Proof}: Let ${\cal B'}_{\sigma ,k}^R $ be the ``signed'' elements 
of ${\cal B}_{\sigma ,k}^R $. We first note that ${\cal S}$ continues
to be a crystal basis for $M'$. Next, the 
equivariance of $\iota_R $ shows that for $v_c \in 
{\cal B'}_{\sigma ,k}^R $, 
\[ \stackrel{\sim}{E_{k}^R } (v_c )=  \iota_R^{-1} (\: \stackrel{\sim }
{e} (\iota_R (v_c ))) \]
This proves
that ${\cal B'}_{\sigma ,k}^R$ is indeed a crystal basis for $V_{\sigma 
,k}^R$. 

Thus, keeping in mind that the signs are alloted by our
global $sign$-function and, by considering all $\sigma$ 
and all $k$, we obtain the assertion. \qed

\eat{Next note that if $\pm v_d =\stackrel{\sim}{E_k^R } (v_c )$ then
$c^* =d^* $. In other words, the action of $\stackrel{\sim}{E_k^R }$ 
stay within ${\cal B'}_{\sigma ,k}^R $ and signs alloted by our global
$sign$-function will never conflict.}

We now define our final global sign $sign^* (b)$ as follows.
Firstly, let $S=\{ b \mid F_i^L v_b =F_k^R v_b =0 \}$. These are 
the lowest weight vectors for both the left and the right action. 
We see that:
\begin{itemize}
\item For any $b \in S$, we have $b^* =b$.
\item If $wt_i (b)$ denotes the cardinality of the set 
$\{ (i,k) | (i,k) \in b \} $, then $wt_1 (b) \leq \ldots \leq wt_m (b)$.
\end{itemize}

We define $sign^*(b)=sign(b)$ for all $b$ such that $b^* \in S$. 
Next, for 
a $c$ such that $c^* \not \in S$, we inductively (by $(wt_i )$ above) 
define $sign^* (c)=sign^* (
\stackrel{\sim}{F_i^L } (v_{c}))$ where  
$\stackrel{\sim}{F_i^L } (v_{c})\neq 0$.
By the commutativity of $\stackrel{\sim}{F_i^L}
$ with $\stackrel{\sim}{F_k^R}$, we see that $sign^* (c)$ is well 
defined over all $K$-subsets of $[m]\times [n]$. 

Let $v^*_b =sign^* (b) \cdot v_b $ and let ${\cal B}^* =\{ v^*_b | v_b \in {\cal B} \}$.

\begin{prop} \label{prop:cb}
The elements ${\cal B}^*$ is a crystal basis 
for $\wedge^K (\C^{mn})$ for the action of both $U_q (gl_n )$ and $U_q
(gl_m )$. In other 
words $\stackrel{\sim }{E_k^R} (v^*_c ) \in {\cal B}^* \cup 0 $.
and $\stackrel{\sim }{E_i^L } (v^*_c ) \in {\cal B}^* \cup 0 $.
\end{prop}

\noindent
{\bf Proof}: The proof follows from the commutativity condition and 
the well-defined-ness of $sign^* $. \qed

\begin{ex} \label{ex:wedge2}

Let us consider $\wedge^2 (\C^{2 \times 2})$ whose six elements, their
matrix notation, and signs are given below:

\[ \begin{array}{ccccc} \hline 
c & matrix & c^* & sign^* (c^* ) & sign^* (c) \\ \hline \hline \\
\vspace{2mm}
\young(1,2) & \young(10,10)& \young(01,01) & 1 & 1 \\ \vspace{2mm}
\young(1,3) & \young(11,00)& \young(11,00) & 1 & 1 \\ \vspace{2mm}
\young(1,4) & \young(10,01)& \young(01,01) & 1 & 1 \\  \vspace{2mm}
\young(2,3) & \young(01,10)& \young(01,01) & 1 & -1 \\  \vspace{2mm}
\young(2,4) & \young(00,11)& \young(00,11) & 1 & 1 \\  \vspace{2mm}
\young(3,4) & \young(01,01)& \young(01,01) & 1 & 1 \\  \\
\hline
\end{array} \]
\end{ex}

For a $b \subseteq  [m]\times [n]$ define the (wedge) {\bf left word}
$WLW(b)$ as the $i$-indices of all elements $(i,k)\in b$, read 
bottom to top within a column, reading the columns left to right. 
Similarly, define the {\bf right 
word} $WRW(b)$ as the $k$-indices of all elements $(i,k) \in b$, 
read right to left within a row, reading the rows from bottom to top. 
For a word $w$, let $rs(w)$ be 
the Robinson-Schenstead tableau associated with $w$, when read from
left to right. Define the {\bf left tableau} $WLT(b)=rs(WLW(b))$ and 
the {\bf right tableau} as $WRT(b)=rs(WRW(b))$.

\begin{ex}
Let $m=3$ and $n=4$ and let $b=\{ 1,3,5,6,9,10,11 \}$. 
\[ \begin{array}{c} 
\young(1001,0101,1110) \: \, \hspace*{1cm}  WLW(b)=3132321 \: \,
\hspace*{1cm}   WRW(b)= 3214241 \\
WLT(b)=\young(112,233,3) \hspace*{1cm} WRT(b)=\young(114,22,34) \\
\end{array}
\]
\end{ex}

For semi-standard tableau, recall the crystal operators
$\stackrel{\sim}{e_i^T} , \stackrel{\sim}{f_i^T }$, see 
for example, \cite{kashiwara}. These crystal operators may be connected 
to our crystal operators via the following proposition. 
This obtains the result in \cite{danil} for the $\wedge$-case.

\begin{prop}
For any $v^*_b \in {\cal B}^* $ the crystal basis for 
$\wedge^K (\C^{mn})$ as above, we have:
\begin{itemize}
\item If $\stackrel{\sim}{E_i^L}(v^*_b )= v^*_c $ then
$\stackrel{\sim}{e_i^T} (WLT(b))=WLT(c)$. 
\item If $\stackrel{\sim}{E_k^R}(v^*_b )= v^*_c $ then
$\stackrel{\sim}{e_k^T} (WRT(b))=WRT(c)$. 
\end{itemize}
A similar assertion holds for the $\stackrel{\sim}{F}$-operators.
\end{prop}
\section{The module $V_{\lambda }$} \label{sec:lambda}

We have thus seen the algebra maps $\phi_L :U_q (gl_m ) \rightarrow 
U_q (gl_{mn}) \rightarrow End_{\C (q)} (\wedge^k (\C^{mn}))$ 
and $\phi_R :U_q (gl_n ) \rightarrow 
End_{\C (q)} (\wedge^k (\C^{mn}) )$.
Since the two actions commute, this converts 
$\wedge^k (\C^{mn})$ into a $U_q (gl_m )\otimes U_q (gl_n )$
-module. Also note that at $q=1$, we have the factorization:
\[ \begin{array}{c}
\phi_L (1) :U_1 (gl_m ) \rightarrow 
U_1 (gl_{mn}) \rightarrow End_{\C} (\wedge^k (\C^{mn})) \\ 
\phi_R (1) :U_1 (gl_n ) \rightarrow 
U_1 (gl_{mn}) \rightarrow End_{\C} (\wedge^k (\C^{mn})) \\ 
\end{array} \]

\begin{prop}
The actions $\phi_L ,\phi_R $ convert $\wedge^k (\C^{mn})$ into 
a $U_q (gl_m )\otimes U_q (gl_n )$ module. Furthermore, at $q=1$ 
this matches the  restriction of the $U_1 (gl_{mn})$ action
on $\wedge^k (\C^{mn})$ to $U_1 (gl_m ) \otimes U_1 (gl_n )$. 
\end{prop}

Since, both $U_q (gl_m )$ and $U_q (gl_n )$ are
Hopf-algebras, we see that if $M,N$ are $U_q (gl_{m})\otimes U_q 
(gl_n )$-modules then so is $M \otimes N$. 
The action of $U_q (gl_m )$ on $M \otimes N$ defined by

\[ \Phi_L : 
U_q (gl_m ) \stackrel{\Delta }{\rightarrow } U_q (gl_m ) \otimes 
U_q (gl_m ) \stackrel{\phi_L \otimes \phi_L }{\rightarrow } U_q (gl_{mn}) \otimes 
U_q (gl_{mn}) \rightarrow End_{\C (q)} (M \otimes N) \]

In the case $M,N$ are $U_q (gl_{mn})$-modules, we also have: 
\[ \Phi'_L : U_q (gl_m ) \stackrel{\phi_L }{\rightarrow } 
U_q (gl_{mn} )  
\stackrel{\Delta }{\rightarrow } U_q (gl_{mn}) \otimes 
U_q (gl_{mn}) \rightarrow End_{\C (q)} (M \otimes N) \]

We may similarly define $\Phi_R $ 
\[ \Phi_R : U_q (gl_n ) \stackrel{\Delta }{\rightarrow } U_q (gl_n) \otimes 
U_q (gl_n ) \stackrel{\phi_R \otimes \phi_R  }{\longrightarrow } End_{\C (q)}
(M \otimes N)  \]
Again, if $M,N$ are $U_q (gl_{mn})$-modules, we have at $q=1$:
\[ \Phi'_R (1): U_1 (gl_n ) \stackrel{\phi_L (1)}{\rightarrow } 
U_1 (gl_{mn} )  
\stackrel{\Delta }{\rightarrow } U_1 (gl_{mn}) \otimes 
U_1 (gl_{mn}) \rightarrow End_{\C } (M \otimes N) \]

\begin{prop} \label{lemma:commute} 
\noindent
\begin{itemize}
\item If $M,N$  are $U_q (gl_{m})\otimes U_q (gl_n )$-modules then 
so is $M\otimes N$, interpreted
as $U_q (gl_m )\otimes U_q (gl_n )$ module through $\Phi_L $
and $\Phi_R $. 
\item The maps $\Phi_L =\Phi'_L $ and $\Phi_R =\Phi'_R $ when $q=1$.  
Thus $\Phi_L $ and $\Phi_R $ are deformations of the action of 
$U_1 (gl_{mn})$ restricted to $U_1 (gl_m ) \otimes U_1 (gl_n )$.
\end{itemize}
\end{prop}

\noindent
The proof of the first part is obvious. For the second part notice
that for $q=1$ both $\phi_L $ and $\phi_R $ match the classical
injections (algebra homomorphisms) of $U_q (gl_m )$ 
(or $U_q (gl_n ))$) into $U_q (gl_{mn})$.

Unless otherwise stated, for $U_q (gl_{mn})$-modules $M,N$, the 
$U_q (gl_m )$ and $U_q (gl_n )$ structure on $M\otimes N$ 
will be that arising from
$\Phi_L $ and $\Phi_R $. \ignore{Let ${\cal M}$ be the collection 
of $U_q (gl_{mn})$-modules treated as $U_q (gl_m )$ and 
$U_q (gl_n )$-modules as above (i.e., via $\Phi_L $ and $\Phi_R $). }

\begin{lemma}
For the module $\wedge^k (\C^{mn}) $ as a $U_q (gl_m )\otimes 
U_q (gl_n )$-module, we have:
\[ \wedge^k (\C^{mn}) =\sum_{\lambda} V_{\lambda} (\C^m ) \otimes V_{\lambda'}
(\C^n ) \]
where $|\lambda |=k$.
\end{lemma}

\noindent
The proof is clear by setting $q=1$. \qed

\noindent
Next, for a $U_1 (gl_{mn})$-module
$V_{\lambda }$ and the standard embedding $U_1 (gl_m )\otimes 
U_1 (gl_n)$, let
\[ V_{\lambda } (\C^{mn})= \oplus_{\alpha ,\beta } \: n^{\lambda }_{\alpha
,\beta } V_{\alpha } (\C^m ) \otimes V_{\beta } (\C^n ) \]

\begin{lemma} \label{lemma:psi}
For $a,b\in \Z$, consider $\wedge^{a+1} (\C^{mn})  \otimes
\wedge^{b-1} (\C^{mn}) $ and 
$\wedge^a (\C^{mn})  \otimes \wedge^b (\C^{mn}) $ as $U_q (gl_m ) \otimes U_q (gl_n
)$-modules. Then there exists an 
$U_q (gl_m )\otimes U_q (gl_n )$-equivariant injection $\psi_{a,b}$:
\[ \psi_{a,b} : \wedge^{a+1} (\C^{mn})  \otimes \wedge^{b-1}
(\C^{mn}) \rightarrow 
\wedge^a (\C^{mn})  \otimes \wedge^b (\C^{mn})  \]
If $\lambda $ is the shape of two columns sized $a$ and $b$ then
the co-kernel $cok(\psi_{a,b})$ may be written as:
\[ cok (\psi_{a,b})=\oplus_{\alpha ,\beta } \: n^{\lambda }_{\alpha ,\beta}
V_{\alpha } (\C^m ) \otimes V_{\beta }(\C^n ) \]
\end{lemma}

\noindent
{\bf Proof}: For $q=1$ the above map is a classical construction (see,
e.g., \cite{fultonrepr}). 
This implies that for general $q$, the multiplicity of the 
$U_q (gl_m )\otimes U_q (gl_n )$-module
$V_{\alpha }(\C^m ) \otimes V_{\beta } (\C^n )$ in $\wedge^{a+1} (\C^{mn}) 
\otimes \wedge^{b-1} (\C^{mn}) $ does not exceed that in
$\wedge^a (\C^{mn}) \otimes 
\wedge^b (\C^{mn}) $. Whence a suitable $\psi_{a,b}$ may be constructed 
respecting  the isotypical components of both modules. The second
assertion now follows. \qed

We now propose a recipe for the construction of the $U_q (gl_m )\otimes 
U_q (gl_n )$ module $W_{\lambda }$. Let $\lambda' =[\mu_1 ,\ldots ,
\mu_r ]$, i.e., $\lambda $ has $r$ columns of length $\mu_1 ,
\ldots ,\mu_r $. Let $C^k $ the the collection
of all columns of size $k$ with strictly increasing entries from
the set $[mn]$. For $a\geq b$ and $c\in C^a $ and $c'\in C^b $, we 
say that $c \leq c'$ if for all $1\leq i \leq b$, we have $c(i)\leq 
c'(i)$. A basis for $W_{\lambda }$ will be the set $SS(\lambda ,mn)$, 
i.e., semi-standard tableau of shape $\lambda $ with entries in $[mn]$.
We interpret this basis as $X^{\lambda } \subseteq 
Z^{\lambda}=\prod_i C^{\mu_i }$. In other words, 
\[ X^{\lambda }=\{ [c_1 ,\ldots ,c_r ] | c_i \in C^{\mu_i }, c_i \leq c_{i+1} \} \]
We call $X^{\lambda }$ as {\bf standard} and 
$Y^{\lambda }=Z^{\lambda} -X^{\lambda} $ as non-standard. We
represent $\wedge^p (\C^{mn}) $ as in \cite{LT}, with the basis $C^p $ and 
construct $M=\otimes_i \wedge^{\mu_i } (\C^{mn}) $ with the basis 
$Z^{\lambda}$. 
Note that $M$ is a $U_q (gl_m ) \otimes U_q (gl_n )$-module. 

Recall the maps:
\[ \psi_{a,b} : \wedge^{a+1} \otimes \wedge^{b-1} \rightarrow \wedge^a \otimes 
\wedge^b \]
Let $Im_{a,b}$ be the image of $\psi_{a,b}$. Define 
\[ \begin{array}{rcl}
{\cal S}_i &=& \wedge^{\mu_1 } \otimes \ldots \otimes \wedge^{\mu_{i-1}} 
\otimes Im_{\mu_i .\mu_{i+1}} \otimes \wedge^{\mu_{i+2}} \otimes \ldots \otimes \wedge^{\mu_r} \\
{\cal S} &=& {\cal S}_1 + \ldots {\cal S}_{r-1} \subseteq M 
\end{array} \]

We call ${\cal S}$ as the {\em straightening laws} for the shape $\lambda $. 
Note that ${\cal S}$ is a $U_q (gl_m )\otimes U_q (gl_n )$-submodule of $M$. 
We may conjecture that a suitable family of $\psi $'s exist so that the desired
module $W_{\lambda }$ is indeed the quotient $M/{\cal S}$ and that the standard
tableau $X^{\lambda}$ form a basis.

\section{The construction of $\psi_{a,b}$} \label{sec:psi}

The structure of the two-column $W_{\lambda}$ and ${\cal S}$ in the general 
case, depend intrinsically on the straightening laws $\psi_{a,b}$ 
(for various $a,b$) of Lemma \ref{lemma:psi}. 
In this section we will construct a family of maps: 
\[ \psi_{a,b} : \wedge^{a+1} \otimes \wedge^{b-1} \rightarrow 
\wedge^a \otimes \wedge^b \]
These maps will have the following important properties:
\begin{itemize}
\item $\psi_{a,b}$ will be $U_q (gl_m )\otimes U_q (gl_n )$-equivariant,
and 
\item at $q=1$, they will also be $U_1 (gl_{mn})$-equivariant and 
will match the standard resolution.
\end{itemize}

This is done in three steps:
\begin{itemize}
\item First, the construction of equivariant maps $\psi_a :\wedge^{a+1}
\rightarrow \wedge^a \otimes \wedge^1 $ and $\psi'_a : \wedge^{a+1} 
\rightarrow \wedge^1 \otimes \wedge^a $.
\item Next, for a module map $\mu :A \rightarrow B$, the construction 
of the ``adjoint'' $\mu^* : B \rightarrow A$.
\item Finally constructing $\psi_{a,b}$ using $\psi_a $ and $\psi^*_b$.
\end{itemize}

We first begin with the adjoint.

\subsection{Normal bases}

Let us fix the basis $B=\{ v_c | c \subseteq [mn] , |c|=k \}$ 
as the basis of $\wedge^k (\C^{mn}) $. We define an inner product on 
$\wedge^k (\C^{mn}) $ as follows. For elements $v_c ,v_{c'} \in 
\wedge^k (\C^{mn}) $,  
let $\langle v_c ,v_{c'} \rangle =\delta_{c,c'}$. In other words,
the inner product is chosen so that $B$ are ortho-normal.

Abusing notation slightly, we denote, for example 
by $\langle E_i^L c,c' \rangle $ as short-form for 
$\langle E_i^L (v_c ) ,v_{c'} \rangle$. We have the {\bf EF-Lemma}:

\begin{lemma} \label{lemma:EF}
For the action of $U_q (gl_m )$ and $U_q (gl_n )$ as above, on 
$\wedge^k (\C^{mn}) $ as above, we have:
\[ q^{-1} q^{h_i^L  }(v_{c'}) \langle E_i^L c,c' \rangle = q q^{h_i^L }
(v_c ) 
\langle E_i^L c,c' \rangle =\langle F_i^L c',c \rangle \] 
\[ q^{-1} q^{h_i^R  }(v_{c'}) \langle E_i^R  c,c'\rangle = q q^{h_i^R }(v_c ) 
\langle E_i^R c,c' \rangle =\langle F_i^R c',c \rangle \] 
\end{lemma}

\noindent
{\bf Proof}: We have:
\[ E^L_i (v_c )=( e_i +q^{-h_i}e_{m+i} + \ldots
    (\prod_{j=0}^{n-2} q^{-h_{jm+i}} ) e_{(n-1)m+i} ) v_c \]
Now, by examining the $gl_{mn}$-weights of $c,c'$, 
exactly one of these terms will lead to $v_{c'}$, and so 
\[ \langle E_i^L c, c'\rangle v_{c'}= 
 (\prod_{j=0}^{k} q^{-h_{jm+i}} ) e_{((k+1)m+i} ) v_c =  
 (\prod_{j=0}^{k} q^{-h_{jm+i}} (v_{c'} )) \cdot v_{c'} \]

Now, we see that:
\[ F^L_i (v_{c'})= ((\prod_{j=1}^{n-1} q^{h_{jm+i}}) f_i +
    \ldots + q^{h_{(n-1)m+i}}f_{(n-2)m+i}+ f_{(n-1)m+i} )v_{c'} \]
It must be the $f_{(k+1)m+i}$ term that led to $v_c $. Whence, we have:
\[ \langle F_i^L c' ,c \rangle 
=((\prod_{j=k+2}^{n-1} q^{h_{jm+i}} ) \cdot f_{(k+1)m+i }) v_{c'} =
(\prod_{j=k+2}^{n-1} q^{h_{jm+i}} ) v_{c} \]
But since $c,c'$ differ only in the entry $(k+1)m+i$, we have
\begin{itemize}
\item $q^{h_{(k+1)m+i}} (v_c ) =q^{-1}$
and $q^{h_{(k+1)m+i}} (v_{c'} ) =q$.
\item 
$ (\prod_{j=0}^{k} q^{-h_{jm+i}} (v_{c'} ))  
= (\prod_{j=0}^{k} q^{-h_{jm+i}} (v_{c} )) $
\item $q^{h_i^L} (v_c )= \prod_{j=0}^{n-1} q^{h_{jm+i}} v_c $
\end{itemize}

Finally,
\[ \begin{array}{rcl}
q q^{h_i^L } (v_c ) \langle E_i^L c,c' \rangle  &=& 
q \prod_{j=k+1}^{n-1} q^{h_{jm+i}} (v_c ) \\
&=& 
q q^{h_{(k+1)m+i}} (v_c ) \prod_{j=k+2}^{n-1} q^{h_{jm+i}} (v_c ) \\
&=& \langle F_i^L c',c \rangle \\ 
\end{array} \]
Other assertions are similarly proved. \qed

\begin{defn}
Let $A$ be a $U_q (gl_m ) \otimes U_q (gl_n )$ module, and let
${\cal A}=\{ a_1 ,\ldots ,a_r \} $ be a basis of $A$ of weight 
vectors. Define an inner product 
$\langle \cdot , \cdot \rangle $ on $A$ making ${\cal A}$ orthogonal. 
We say that ${\cal A}$ is {\em normal} if the EF-lemma Lemma 
\ref{lemma:EF} holds (with $a,a'\in {\cal A}$ replacing $c,c'$).
\end{defn}

\begin{lemma} \label{lemma:normaltensor}
Let $A,B$ be $U_q (gl_m ) \otimes U_q (gl_n )$-modules such
that ${\cal A}=\{ a_1 ,\ldots ,a_r \}$ and ${\cal B}= \{ b_1 ,
\ldots , b_s \}$ are normal bases for $A$ and $B$ respectively.
Then ${\cal A} \otimes {\cal B}$ is a normal basis for $A \otimes 
B$ with the inner product $\langle a\otimes b ,a'\otimes b' \rangle 
=\delta_{a \otimes b, a'\otimes b' }$.
\end{lemma}

\noindent
{\bf Proof}: Let consider the element $a\otimes b$, and the 
elements $a'\otimes b$ and $a\otimes b'$ such that $a'$ appears 
in $E_i^L a $ and $b'$ appears in $E_i^L b$. 

We see that:
\[ \begin{array}{rcl}
q\cdot q^{h_i^L } (a \otimes b) \langle E_i^L (a \otimes b), a' \otimes b \rangle 
&=& q\cdot q^{h_i^L } (a \otimes b) \langle (E_i^L \otimes 1 +q^{-h_i^L } \otimes E_i^L ) 
(a \otimes b) ,a' \otimes b \rangle \\
&=& q\cdot q^{h_i^L } (a \otimes b) \langle (E_i^L \otimes 1 ) 
(a \otimes b) ,a' \otimes b \rangle \\
&=& q\cdot q^{h_i^L } (a \otimes b) \langle E_i^L a ,a' \rangle \\
&=& q^{h_i^L } (b) [q\cdot q^{h_i^L } (a)  \langle E_i^L a ,a' 
\rangle ] \\
&=& q^{h_i^L } (b) \langle F_i^L a' ,a \rangle \\
\end{array} \]
On the other hand, we have:
\[ \begin{array}{rcl}
\langle F_i^L (a'\otimes b) ,a \otimes b \rangle 
&=& \langle (F_i^L \otimes q^{h_i^L } + 1 \otimes F_i^L )
(a'\otimes b), a \otimes b \rangle \\
&=& \langle (F_i^L \otimes q^{h_i^L } )
(a'\otimes b), a \otimes b \rangle \\
&=& q^{h_i^L } (b) \langle F_i^L a' ,a \rangle \\
\end{array} \]

\noindent
Other cases are similar. \qed

Let $\Xi $ be the $\Z$-submodule generated by $\epsilon_i^L $ and
$\epsilon_j^R $.
Let $\chi $ be a $\Xi $-weight and let $\chi' =\chi +h_i^L $. For a 
module $A$ with a normal base ${\cal A}$, let $A_{\chi }$ be the 
weight-space of weight $\chi $.
We see that $E_i^L : A_{\chi } \rightarrow A_{\chi'}$ , while 
$F_i^L : A_{\chi '} \rightarrow A_{\chi }$. Let $a_{\chi}$ be the 
column-vector of elements of ${\cal A}$ of weight $\chi $. Let us 
define matrices $E^A,F^A$ as:
\[ E^A a_{\chi '}=E_i^L a_{\chi} \mbox{ \hspace*{2cm}} F^A a_{\chi}=
F_i^L a_{\chi'} \]

\noindent
By the EF-lemma (i.e., Lemma \ref{lemma:EF}), 
\[ q\cdot q^{<\chi .h_i^L >} E^A =(F^A)^T \]

Now, let $A$ and $B$ be $U_q (gl_m )\otimes U_q (gl_n )$ with normal 
bases ${\cal A}$ and ${\cal B}$ respectively. Let $\mu: A \rightarrow
B$ be an equivariant map and let $\mu_{\chi}$ be a matrix such that:
\[ \mu a_{\chi} =\mu_{\chi } b_{\chi} \]
Equivariance implies:
\[ \mu \cdot E_i^L a_{\chi} = \mu \cdot E^A a_{\chi'}
= E^A \mu_{\chi } b_{\chi'} \]
\[ E_i^L \cdot \mu a_{\chi} = E_i^L \mu_{\chi} b_{\chi} =
\mu_{\chi'} E^B b_{\chi'} \]
\noindent
Or in other words,
\[ E^A \mu_{\chi}  = \mu_{\chi'} E^B  \mbox{\hspace*{2cm}}
F^A \mu_{\chi'}  = \mu_{\chi} F^B \]
Transposing the second equivariance condition, we get:
\[ (F^A \mu_{\chi'})^T  = (\mu_{\chi} F^B )^T \] 
We may simplify this as:
\[ \mu_{\chi'}^T (F_A)^T =(F_B)^T \mu_{\chi}^T \]
and further:
\[  q\cdot q^{< \chi . h_i^L >} \mu_{\chi'}^T E_A = q\cdot q^{< \chi .
h_i^L >} E_B \mu_{\chi}^T \]
i.e., finally:
\[  \mu_{\chi'}^T E_A = E_B \mu_{\chi}^T \]
We may similarly prove that
\[  \mu_{\chi}^T F_A = F_B \mu_{\chi'}^T \]
Both these observations immediately imply:
\begin{prop} \label{prop:transpose}
Let $\mu :A \rightarrow B$ be an equivariant map, and let $\mu_{\chi}$ be
defined as above. We construct the map $\mu^* :B \rightarrow A$ as 
follows. Define $\mu^* $ such that:
\[ \mu^* b_{\chi} =\mu_{\chi}^T a_{\chi } \]
Then $\mu^* :B \rightarrow A$ is equivariant.
\end{prop}

\subsection{The Construction of $\psi_a$} 
In this section we construct the $U_q(gl_m) \otimes
U_q(gl_n)$-equivariant maps
\begin{eqnarray*} 
\psi_{a} : \wedge^{a+1}  \rightarrow \wedge^a \otimes \wedge^1 \\
\psi'_{a} : \wedge^{a+1}  \rightarrow \wedge^1 \otimes \wedge^a 
\end{eqnarray*} 

Note that $\wedge^1 = \C^{mn} = \C^m \otimes \C^n$. For convenience,
we identify $[mn]$ with $[m]\times [n]$. Under this identification,
an element $(i,j) \in [m]\times [n]$ maps to the element
$m*(j-1)+i$.  

In this notation, the natural basis for the representation 
$\wedge^k = \wedge^k (\C^{mn}) $ is parametrized by subsets of 
$[m]\times [n]$ with $k$ elements.

Recall that, as a $U_q (gl_m) \otimes U_q (gl_n)$-module, we have
\[ \wedge^k (\C^{mn}) =\sum_{\lambda} V_{\lambda} (\C^m ) \otimes V_{\lambda'}
(\C^n ) \]
where $|\lambda |=k$. Further, $\lambda$ has atmost $m$ parts and
$\lambda'$ has atmost $n$ parts, that is, the shape $\lambda$ fits
inside the $m\times n$ `rectangle'.

For a shape $\lambda = (\lambda_1, \ldots, \lambda_m)$ with
$\lambda' = (\lambda'_1, \ldots, \lambda'_n)$, consider 
the subset $c_\lambda \subset [mn]$
defined as:
\[ 
\begin{array}{lllll}
c_\lambda & = & \{ & 1, m+1,  \ldots \ldots \ldots, m*(\lambda_1-1)+1, & \\
          &   &    & 2, m+2, \ldots \ldots, m*(\lambda_2-1)+2, & \\
	  &   &    & \ldots \ldots & \\
	  &   &    &  m, 2m, \ldots, m*(\lambda_m-1) & \}
\end{array}
\]
Equivalently, 
\[
\begin{array}{lllll}
c_\lambda & = & \{ & 1, 2,  \ldots \ldots \ldots, \lambda'_1, & \\
          &   &    & m+1, m+2, \ldots \ldots, m+\lambda'_2, & \\
	  &   &    & \ldots \ldots & \\
	  &   &    & m*(n-1)+1, m*(n-1)+2, \ldots, m*(n-1)+\lambda'_n & \}
\end{array}
\]
Under the identification of $[mn]$ with $[m]\times[n]$, we have
\[
c_\lambda = \left\{ 
(i,j) \mid 1\leq i \leq \lambda'_j, 1\leq j \leq \lambda_i \right\}
\]
We slightly abuse the notation and write $(i,j) \in \lambda$ as 
a short-form for $(i,j) \in c_\lambda$.

With this notation, we have the following important lemma:
\begin{lemma}
Consider the $U_q(gl_m)\otimes U_q(gl_n)$-module $\wedge^k (\C^{mn})$. 
For a shape $\lambda$ which fits in the $m \times n$ rectangle
with $|\lambda|=k$, the weight vector $v_{c_\lambda} \in \wedge^k$ 
is the highest $U_q(gl_m)\otimes U_q(gl_n)$-weight vector of weight
$(\lambda, \lambda')$.
\end{lemma}
\proof
The lemma follows from the observation that
$E_i^L(v_{c_\lambda}) = E_j^R (v_{c_\lambda}) = 0$ for all $i,j$.
\qed

Now we turn our attention to the construction of the 
$U_q(gl_m)\otimes U_q(gl_n)$-equivariant  map 
$$\psi_{a} : \wedge^{a+1}  \rightarrow \wedge^a \otimes \wedge^1 $$

As a $U_q(gl_m) \otimes U_q(gl_n)$-module, we have the following
decomposition
\[ \wedge^{a+1}=
\sum_{\lambda : |\lambda|=a+1} V_{\lambda} (\C^m ) \otimes 
V_{\lambda'} (\C^n ) \]
Moreover, $v_{c_\lambda}$ is the highest-weight vector for the
$U_q(gl_m)\otimes U_q(gl_n)$-submodule 
$V_{\lambda} (\C^m ) \otimes V_{\lambda'} (\C^n)$ of $\wedge^{a+1}$.

Thus, in order to construct the $U_q(gl_m)\otimes
U_q(gl_n)$-equivariant map $\psi_{a}$, we need to simply
define the images $\psi_a(v_{c_{\lambda}})$ inside 
$\wedge^a \otimes \wedge^1$. Moreover the vector 
$\psi_a(v_{c_{\lambda}})$ should be a highest-weight vector
of weight $(\lambda, \lambda')$. Note that, unlike $\wedge^{a+1}$,
$\wedge^a \otimes \wedge^1$ is not multiplicity-free. 
Below, we outline the construction of a highest-weight vector 
(upto scalar multiple) $v_\lambda$
of weight $(\lambda, \lambda')$ inside $\wedge^a \otimes \wedge^1$.

We begin with some notation. As before, fix a shape $\lambda$
which fits in the $m \times n$ rectangle with $|\lambda|=a+1$.
Write $\lambda = (\lambda_1, \ldots, \lambda_m)$ with
$\lambda' = (\lambda'_1, \ldots, \lambda'_n)$ and 
\[
c_\lambda = \left\{ 
(i,j) \mid 1\leq i \leq \lambda'_j, 1\leq j \leq \lambda_i \right\}
\]

For $(i,j) \in \lambda$, we set 
\begin{eqnarray*}
t_{i,j} = v_{c_\lambda - \{i,j\}} \in \wedge^a \\
\chi_{i,j}=v_{\{(i,j)\}} \in \wedge^1
\end{eqnarray*}
In other words, $t_{i,j}$ is the vector in $\wedge^a$
corresponding to the subset obtained from the subset $c_{\lambda}$
by removing the element $(i,j) \in \lambda$. Further, 
$\chi_{i,j}$ is the vector in $\wedge^1$ corresponding to
the singleton set containing the element $(i,j)$. Below, we
abuse notations and denote by $t_{i,j}$ and $\chi_{i,j}$ also
the subsets that correspond to these vectors.

\begin{lemma} For $(i,j) \in \lambda$, $1 \leq k < m$, $1 \leq l < n$,
\begin{itemize}
\item $E^L_k (t_{i,j}) = 0$ if $i \neq k$.
\item $E^L_i (t_{i,j}) = t_{i+1,j}$ if $(i+1,j) \in \lambda$ and 0
otherwise.
\item If $(i+1,j) \in \lambda$, 
$q^{-h_i^L} (t_{i+1,j}) = q^{\lambda_{i+1}-\lambda_i-1} t_{i+1,j}$.
\item $E^R_l (t_{i,j}) = 0$ if $j \neq l$.
\item $E^R_j (t_{i,j}) = (-1)^{\lambda'_j-1}q^{\lambda'_{j+1}-\lambda'_j}
t_{i,j+1}$ if $(i,j+1) \in \lambda$ and 0
otherwise.
\item If $(i,j+1) \in \lambda$,
$q^{-h_j^R} (t_{i,j+1}) = q^{\lambda'_{j+1}-\lambda'_j-1} t_{i,j+1}$.
\end{itemize}
\end{lemma}
\proof Let $k \neq i$ and consider 
$E^L_k (t_{i,j})$. Note that, for all $j'$, 
if $(k+1, j') \in t_{i,j}$, then $(k, j') \in t_{i,j}$. Thus,
by definition of $E^L_k$, we have $E^L_k (t_{i,j}) = 0$.

Now consider $E^L_i (t_{i,j})$.  Note that $(i,j) \not\in t_{i,j}$.
If $(i+1, j) \in \lambda$, then $(i+1,j) \in t_{i,j}$. Further,
for all $j' < j$, if $(i+1,j') \in t_{i,j}$ then $(i, j') \in t_{i,j}$.
Thus, by definition $E^L_i(t_{i,j})$ operates only at the
position $(i+1, j)$ if $(i+1, j) \in \lambda$ and produces
the subset $t_{i+1,j}$.

Now we assume that $(i+1,j) \in \lambda$, and evaluate
$q^{-h_i^L} (t_{i+1,j})$. Note that, except for $(i+1,j)$, 
$(i+1,j') \in t_{i+1,j}$  for $1 \leq j' \leq \lambda_{i+1}$.
Also, for $j' > \lambda_{i+1}$, $(i+1, j') \notin t_{i+1,j}$. Thus 
$q^{\epsilon^L_{i+1}} (t_{i+1,j}) = q^{\lambda_{i+1}-1}$. Similarly,
$q^{\epsilon^L_{i}} (t_{i+1,j}) = q^{\lambda_{i}}$. Therefore,
\[q^{-h_i^L} (t_{i+1,j}) = q^{\lambda_{i+1}-\lambda_i-1} t_{i+1,j}\]

It is easy to that $E^R_l (t_{i,j}) = 0$ if $j \neq l$. So, we
turn our attention to $E^R_j (t_{i,j})$. Note that, 
for $i'$ such that 
$\lambda'_{j+1} < i' \leq \lambda'_j$, $(i',j) \in t_{i,j}$ and
$(i',j+1) \not\in t_{i,j}$. For other values of $i'$ except $i$, 
either both or none of $(i',j)$ and $(i',j+1)$ belong to $t_{i,j}$.
Therefore, as expected,
$E^R_j(t_{i,j})$ operates only at the position $(i,j+1)$ if
$(i,j+1) \in \lambda$. Further, by definition of $E^R_j$, if
$(i,j+1) \in \lambda$, we have 
$$E^R_j (t_{i,j}) = (-1)^{\lambda'_j-1}q^{\lambda'_{j+1}-\lambda'_j}t_{i,j+1} \mbox{~if~} (i,j+1) \in t_{i,j}$$
The sign $(-1)^{\lambda'_j-1}$ results from the fact that exactly
$\lambda'_j-1$ elements of $[mn]$ strictly in the range from $(i,j)$ to
$(i,j+1)$ belong to $t_{i,j}$.

We skip the proof for the last assertion as it follows from a similar
reasoning applied earlier for the left $E$-operator.
\qed

\begin{lemma} For $(i,j) \in \lambda$,
\begin{itemize}
\item \[
E^L_i (t_{i,j} \otimes \chi_{i,j}) =  
\left\{
  \begin{array}{lr} 
  t_{i+1,j}\otimes \chi_{i,j} & \mbox{~if~} (i+1,j) \in \lambda \\
  0 & \mbox{~otherwise~} 
  \end{array}
\right. \\
\]
\item If $(i+1, j) \in \lambda$, then 
\[
E^L_i (t_{i+1,j} \otimes \chi_{i+1,j}) =  
  q^{\lambda_{i+1}-\lambda_i-1} t_{i+1,j}\otimes \chi_{i,j} 
\]
\item 
\[
E^R_j (t_{i,j} \otimes \chi_{i,j}) =  
\left\{
  \begin{array}{lr} 
  (-1)^{\lambda'_j - 1}
  q^{\lambda'_{j+1}-\lambda'_j} t_{i,j+1}\otimes \chi_{i,j} & 
  \mbox{~if~} (i,j+1) \in \lambda \\
  0 & \mbox{~otherwise~} 
  \end{array}
\right.
\]
\item 
If $(i, j+1) \in \lambda$, then 
\[
E^R_j (t_{i,j+1} \otimes \chi_{i,j+1}) =  
  q^{\lambda'_{j+1}-\lambda'_j -1} t_{i,j+1}\otimes \chi_{i,j}
\]
\item For remaining $1 \leq k < m$ and $1 \leq l < n$, 
$E^L_k (t_{i,j} \otimes \chi_{i,j})=
E^R_l (t_{i,j} \otimes \chi_{i,j})= 0$.
\end{itemize}
\end{lemma}
\proof For the first assertion, consider 
\[
E^L_i (t_{i,j} \otimes \chi_{i,j})  =    
   E^L_i (t_{i,j}) \otimes \chi_{i,j} +
   q^{-h_i^L} (t_{i,j}) \otimes E^L_i(\chi_{i,j}) 
\]
As $(i+1,j) \not\in \chi_{i,j}$, 
$E^L_i(\chi_{i,j}) = 0$. Therefore, the claim follows from
the previous lemma.

For the second assertion, let us assume that $(i+1,j) \in \lambda$. Then
\[
E^L_i (t_{i+1,j} \otimes \chi_{i+1,j})  =    
   E^L_i (t_{i+1,j}) \otimes \chi_{i+1,j} +
   q^{-h_i^L} (t_{i+1,j}) \otimes E^L_i(\chi_{i+1,j}) 
\]
Note that, from the previous lemma $E^L_i (t_{i+1,j}) = 0$. Also, 
$E^L_i(\chi_{i+1,j}) = \chi_{i,j}$. Again, using the previous lemma,
we have 
\[
E^L_i (t_{i+1,j} \otimes \chi_{i+1,j})  =    
 q^{\lambda_{i+1}-\lambda_i-1} t_{i+1,j}   \otimes \chi_{i,j} 
\]
The third and fourth assertions are proved in a similar fashion.
\qed

\begin{lemma} Let $v_\lambda \in \wedge^a \otimes \wedge^1$ be
defined as follows:
$$v_\lambda = \sum_{(k,l) \in \lambda} \alpha_{k,l} t_{k,l}\otimes \chi_{k,l}$$
where
$$ \alpha_{k,l}=(-1)^{\lambda'_1 + \ldots + \lambda'_{l-1}+k}
                 q^{k+l-\lambda_k}$$
Then $v_\lambda$ is a highest-weight vector  of weight $(\lambda,
\lambda')$.
\end{lemma}
\proof It is clear that $v_\lambda$ is a weight vector of weight 
$(\lambda, \lambda')$. Below, we show that it is a highest-weight
vector by checking that $E^L_i(v_\lambda) = E^R_j(v_\lambda)= 0$ for
all $i,j$.

Towards this, by previous lemma, we have
\[
\begin{array}{rll}
E^L_i(v_\lambda) = & \sum_{(k,l) \in \lambda} \alpha_{k,l} 
                 E^L_i\left(t_{k,l}\otimes \chi_{k,l}\right) \\
  = & \sum_{(i,l) \in \lambda} \alpha_{i,l} 
                 E^L_i\left(t_{i,l}\otimes \chi_{i,l}\right) +
                 \sum_{(i+1,l) \in \lambda} \alpha_{i+1,l} 
                 E^L_i\left(t_{i+1,l}\otimes \chi_{i+1,l}\right) \\
  = & \sum_{l: (i,l) \& (i+1,l) \in \lambda} 
      \left(  \alpha_{i,l}E^L_i\left(t_{i,l}\otimes \chi_{i,l}\right) +
              \alpha_{i+1,l} 
                 E^L_i\left(t_{i+1,l}\otimes \chi_{i+1,l}\right)
      \right)
\end{array}
\]
For $l$ such that both $(i,l)$ and $(i+1,l)$  are in $\lambda$,
from previous lemma, we have
\begin{eqnarray*}
E^L_i\left(t_{i,l}\otimes \chi_{i,l}\right) =
t_{i+1,l}\otimes \chi_{i,l} \\
E^L_i\left(t_{i+1,l}\otimes \chi_{i+1,l}\right)
=  q^{\lambda_{i+1}-\lambda_i-1} t_{i+1,l}\otimes \chi_{i,l} 
\end{eqnarray*}
Therefore, the coefficent of $t_{i+1,l}\otimes \chi_{i,l}$ in 
$E^L_i(v_\lambda)$ is 
\[
\begin{array}{ll}
= & \alpha_{i,l} + q^{\lambda_{i+1}-\lambda_i-1} \alpha_{i+1,l}\\
= & (-1)^{\lambda'_1 + \ldots + \lambda'_{l-1}+i} q^{i+l-\lambda_i}
    + q^{\lambda_{i+1}-\lambda_i-1}
    (-1)^{\lambda'_1 + \ldots + \lambda'_{l-1}+i+1}
    q^{i+1+l-\lambda_{i+1}}\\
= & (-1)^{\lambda'_1 + \ldots + \lambda'_{l-1}+i} (q^{i+l-\lambda_i}
    - q^{i+l-\lambda_i})\\
= & 0
\end{array}
\]
Thus, $E^L_i(v_\lambda)=0$.
A similar analysis shows that, the coefficient of 
$t_{k, j+1}\otimes \chi_{k,j}$ in $E^R_j (v_\lambda)$ is 
\[
\begin{array}{ll}
= & \alpha_{k,j} (-1)^{\lambda'_j - 1}
  q^{\lambda'_{j+1}-\lambda'_j}+ 
  \alpha_{k,j+1}q^{\lambda'_{j+1}-\lambda'_j -1} 
  \\
= & (-1)^{\lambda'_1 + \ldots + \lambda'_{j-1}+k} q^{k+j-\lambda_k}
    (-1)^{\lambda'_j - 1} q^{\lambda'_{j+1}-\lambda'_j}
    + (-1)^{\lambda'_1 + \ldots + \lambda'_{j}+k} q^{k+j+1-\lambda_{k}}
    q^{\lambda'_{j+1}-\lambda'_j -1} 
    \\
= & 0
\end{array}
\]
This shows that $E^R_j(v_\lambda)=0$ and hence establishes the claim
that $v_\lambda$ is a highest-weight vector in $\wedge^a\otimes \wedge^1$.
\qed

We remark that in the above expression for $v_\lambda$, the coefficient,
$\alpha_{1,1}$, of the term $t_{1,1}\otimes \chi_{1,1}$ has the least
$q$-degree. We may normalize $v_\lambda$ so as to ensure that 
$\alpha_{1,1}=1$ and all the other terms have strictly positive
$q$-degree. 

Now we are ready to define the $U_q(gl_m) \otimes U_q(gl_n)$-
equivariant map
$$\psi_{a} : \wedge^{a+1}  \rightarrow \wedge^a \otimes \wedge^1$$
This is done by simply setting $\psi_a (v_{c_\lambda}) = v_\lambda$. 
It is easily seen that there is a unique 
$U_q(gl_m)\otimes U_q(gl_n)$-equivariant extension of 
$\psi_a$ to all of $\wedge^{a+1}$. Moreover, this extension matches
the classical $U_1(gl_{mn})$-equivariant construction 
at $q=1$.

Next, we prepare towards the construction of the
$U_q(gl_m)\otimes U_q(gl_n)$-equivariant map $\psi'_a$.

\begin{lemma} For $(i,j) \in \lambda$,
\begin{itemize}
\item \[
E^L_i (\chi_{i,j} \otimes t_{i,j}) =  
\left\{
  \begin{array}{lr} 
  q^{-1}
  \chi_{i,j}\otimes t_{i+1,j} & \mbox{~if~} (i+1,j) \in \lambda \\
  0 & \mbox{~otherwise~} 
  \end{array}
\right. \\
\]
\item If $(i+1, j) \in \lambda$, then 
$E^L_i (\chi_{i+1,j}\otimes t_{i+1,j}) =  \chi_{i,j} \otimes t_{i+1,j}$
\item 
\[
E^R_j (\chi_{i,j}\otimes t_{i,j} ) =  
\left\{
  \begin{array}{lr} 
  (-1)^{\lambda'_j - 1}
  q^{\lambda'_{j+1}-\lambda'_j-1} \chi_{i,j}\otimes t_{i,j+1} & 
  \mbox{~if~} (i,j+1) \in \lambda \\
  0 & \mbox{~otherwise~} 
  \end{array}
\right.
\]
\item 
If $(i, j+1) \in \lambda$, then 
\[
E^R_j ( \chi_{i,j+1}\otimes t_{i,j+1}) =  
  \chi_{i,j}\otimes t_{i,j+1}
\]
\item For remaining $1 \leq k < m$ and $1 \leq l < n$, 
$E^L_k (\chi_{i,j} \otimes t_{i,j})=
E^R_l (\chi_{i,j} \otimes t_{i,j})= 0$.
\end{itemize}
\end{lemma}
\proof For the first assertion, consider 
\[
E^L_i (\chi_{i,j}\otimes t_{i,j} )  =    
   E^L_i (\chi_{i,j}) \otimes t_{i,j} +
   q^{-h_i^L} (\chi_{i,j}) \otimes E^L_i(t_{i,j}) 
\]
As $(i+1,j) \not\in \chi_{i,j}$, 
$E^L_i(\chi_{i,j}) = 0$. Further, 
$q^{-h_i^L} (\chi_{i,j}) = q^{-1} \chi_{i,j}$.
Therefore, the claim follows.

For the second assertion, let us assume that $(i+1,j) \in \lambda$. Then
\[
E^L_i (\chi_{i+1,j} \otimes t_{i+1,j})  =    
   E^L_i (\chi_{i+1,j} ) \otimes t_{i+1,j}+
   q^{-h_i^L} (\chi_{i+1,j}) \otimes E^L_i(t_{i+1,j}) 
\]
Note that, $E^L_i (t_{i+1,j}) = 0$. Also, 
$E^L_i(\chi_{i+1,j}) = \chi_{i,j}$. Therefore, 
we have 
\[
E^L_i (\chi_{i+1,j}\otimes t_{i+1,j}  )  =    
 \chi_{i,j}\otimes  t_{i+1,j}   
\]
For the third assertion, consider
\[
\begin{array}{lll}
E^R_j (\chi_{i,j}\otimes t_{i,j} ) &  =    &
   E^R_j (\chi_{i,j}) \otimes t_{i,j} +
   q^{-h_j^R} (\chi_{i,j}) \otimes E^R_j(t_{i,j}) \\
   & = & 
   q^{-1} \chi_{i,j} \otimes E^R_j(t_{i,j})
\end{array}
\]
Recall that, we have 
$$E^R_j (t_{i,j}) = (-1)^{\lambda'_j-1}q^{\lambda'_{j+1}-\lambda'_j}
t_{i,j+1} \mbox{~if~} (i,j+1) \in \lambda \mbox{~and~} 0 \mbox{~otherwise~}$$
Therefore, the claim follows.

For the fourth claim, we assume $(i, j+1) \in \lambda$. Then 
\[
\begin{array}{lll}
E^R_j ( \chi_{i,j+1}\otimes t_{i,j+1}) & =  & 
  E^R_j(\chi_{i,j+1})\otimes t_{i,j+1} +
  q^{-h_j^R}(\chi_{i,j+1})\otimes E^R_j (t_{i,j+1}) \\
  & = & \chi_{i,j}\otimes t_{i,j+1}
\end{array}
\]
The last claim can be easily proved.
\qed

\begin{lemma} Let $v_\lambda \in \wedge^1 \otimes \wedge^a$ be
defined as follows:
$$v_\lambda = \sum_{(k,l) \in \lambda} \beta_{k,l} 
\chi_{k,l} \otimes t_{k,l}$$
where
$$ \beta_{k,l}=(-1)^{\lambda'_1 + \ldots + \lambda'_{l-1}+k}
                 q^{\lambda'_l-k-l}$$
Then $v_\lambda$ is a highest-weight vector  of weight $(\lambda,
\lambda')$.
\end{lemma}
\proof
Clearly, $v_{\lambda}$ is a weight-vector of weight $(\lambda, \lambda')$.
We now check that $E^L_i(v_\lambda)=0$ for all $i$. As expected, this
finally reduces to checking if the following expression, coefficient of
$\chi_{i,l}\otimes t_{i+1, l}$ in $E^L_i(v_\lambda)$, is zero. Towards
this, consider
\[
\begin{array}{ll}
= & q^{-1}\beta_{i,l} + \beta_{i+1,l}\\
= & q^{-1}(-1)^{\lambda'_1 + \ldots + \lambda'_{l-1}+i} q^{\lambda'_l-i-l}
    + (-1)^{\lambda'_1 + \ldots + \lambda'_{l-1}+i+1}
    q^{\lambda'_l-i-1-l}\\
= & 0
\end{array}
\]

Similarly, to check if $E^R_j(v_\lambda)=0$, we need to check if
the following expression, coefficient of
$\chi_{k,j}\otimes t_{k, j+1}$ in $E^R_j(v_\lambda)$, is zero. Towards
this, consider
\[
\begin{array}{ll}
= & \beta_{k,j} (-1)^{\lambda'_j - 1}
  q^{\lambda'_{j+1}-\lambda'_j-1}+ 
  \beta_{k,j+1} \\
= & (-1)^{\lambda'_1 + \ldots + \lambda'_{j-1}+k} q^{\lambda'_j-k-j}
    (-1)^{\lambda'_j - 1} q^{\lambda'_{j+1}-\lambda'_j-1}
    + (-1)^{\lambda'_1 + \ldots + \lambda'_{j}+k} 
    q^{\lambda'_{j+1}-k-j-1}
    \\
= & 0
\end{array}
\]
Thus, we have verified that $E^L_i(v_\lambda)=E^R_j(v_\lambda)=0$ for
all $i,j$. This shows that $v_\lambda$ is a highest-weight vector.
\qed

Now we are ready to define the $U_q(gl_m) \otimes U_q(gl_n)$-
equivariant map
$$\psi'_{a} : \wedge^{a+1}  \rightarrow \wedge^1 \otimes \wedge^a$$
As expected, this is done by simply setting 
$\psi'_a (v_{c_\lambda}) = v'_\lambda$ and taking the
unique $U_q(gl_m)\otimes U_q(gl_n)$-equivariant extension.
Also, as before, this extension matches
the classical $U_1(gl_{mn})$-equivariant construction 
at $q=1$.

Note that $\wedge^{a+1}$ and $\wedge^1 \otimes \wedge^a $ have normal
bases. Whence, by Prop. \ref{prop:transpose}, there is the 
$U_q (gl_m )\otimes U_q (gl_n )$-equivariant map:
$$\psi'^*_{a} : \wedge^1 \otimes \wedge^a \rightarrow  \wedge^{a+1}$$

Finally, we construct $\psi_{a,b}$ as follows:
$$\psi_{a,b}: \wedge^{a+1}\otimes \wedge^{b-1} 
\stackrel{\psi_a \otimes I_{\wedge^{b-1}}}{\longrightarrow}
\wedge^a \otimes \wedge^1 \otimes \wedge^{b-1}
\stackrel{I_{\wedge^{a}} \otimes \psi'^*_{b-1}}{\longrightarrow}
\wedge^{a} \otimes \wedge^{b}$$.

\section{The $Sym^k$ modules} \label{sec:sym}

In this section, we develop the structure of $Sym^k (\C^{mn})$, i.e., 
when $\lambda =k$ or $\lambda '=[1,1,\ldots ,1]$. We start with the module
$M_k =\wedge^1 \otimes \ldots \otimes \wedge^1 $. The straightening laws
${\cal S}$ are generated by exactly the images of $\psi_2 : \wedge^2 
\rightarrow \wedge^1 \otimes \wedge^1 $. 

We use the variables $x_r $ for 
$r=1,\ldots ,mn$ as a basis for $\C^{mn}$ or alternately treat $X$ as a matrix
with the basis $z_{i,j}$ for 
$i=1,\ldots ,m$ and $j=1,\ldots ,n$ with the understanding that 
$z_{i,j}=x_{(j-1)m+i)}$. 
We use $(i,j)$ to mean the element $(j-1)m+i \in [mn]$. We have
an order $(i,j) \leq (i',j')$ which comes from their being elements of $[mn]$.

From the previous section, we 
see that for general $U_q (gl_m ) \otimes U_q (gl_n )$-action, 
$\wedge^2 $ has two highest weight sets $c_{(2)}=
\{ 1,m+1 \}$ and $c_{(1,1)} =\{ 1,2 \}$, and their corresponding vectors 
$v_{c_{(2)}}$ and $v_{c_{(1,1)}}$. 
We see that:
\[ \begin{array}{rcl}
\psi_a (v_{c_{(2)}}) &=&  q \cdot x_1\otimes x_{m+1} - x_{m+1} \otimes x_1 \\
\psi_a (v_{c_{(1,1)}}) &=&  q \cdot x_1\otimes x_2 - x_2 \otimes x_1 \\
\end{array} \]

In the $z$-notation, we thus have:
\[ \begin{array}{l}
q \cdot z_{1,1} \otimes z_{1,2} -z_{1,2} \otimes z_{1,1} \in {\cal S} \\
q \cdot z_{1,1} \otimes z_{2,1} -z_{2,1} \otimes z_{1,1} \in {\cal S} \\
\end{array} \]

The action of $U_q (gl_m )$ and $U_q (gl_n )$ yield more straightening laws
by which {\em non-standard} tuples may be expressed in terms of 
standard tuples. The exact expressions appear below (after dropping the 
$\otimes $): 
\[ \begin{array}{rcl}
z_{i,j+r} z_{i,j}&=& q \cdot z_{i,j} z_{i,j+r} \\
z_{i+r,j} z_{i,j}&=& q \cdot z_{i,j} z_{i+r,j} \\
z_{i,j+s} z_{i+r,j}&=& z_{i+r,j} z_{i,j+s} \\
z_{i+r,j+s} z_{i,j}&=& z_{i,j} z_{i+r,j+s}+ (q-1/q)\cdot z_{i+r,j} z_{i,j+s}\\
\end{array} \]

Denote by $S_k (X)$ the module $M_k /{\cal S}$ and $S(X)$ as $\oplus_k S_k 
(X)$. It has been shown 
elsewhere (see, e.g., \cite{klimyk}), that $S(X)$ is in fact, an associative 
algebra, the so called {\bf quantum matrix algebra} $M_q (X)$. The degree 
$d$-component is indeed exactly spanned by:
\[ z_{1,1}^{d_{1,1}} z_{2,1}^{d_{2,1}} \ldots z_{m-1,n}^{d_{m-1,n}} 
z_{m,n}^{d_{m,n}} \: \: \: \mbox{     where    }\: \: \:  \sum_{i,j} d_{i,j}=d \]
Thus, we have constructed the module $\wedge^k (\C^{mn})$ as a 
$U_q (gl_m ) \otimes U_q (gl_n )$-module. We will now construct a 
crystal basis. 

\noindent
Note that 
\[ \begin{array}{rcl}
  E_j^R (z_{i',j'})&=& \left\{ \begin{array}{lr} 
            z_{i',j} & \: \: \mbox{ if } j'=j+1 \\
            0   & \mbox{otherwise} \end{array} \right. \\ 
&& \\
  F_j^R (z_{i',j'})&=& \left\{ \begin{array}{lr} 
            z_{i',j+1} & \: \: \mbox{ if } j'=j \\
            0   & \mbox{otherwise} \end{array} \right. \\ 
\end{array} \]
We will use the ``standard'' Hopf defined below:
\[ \begin{array}{rcl}
\Delta (e) &=& e \otimes 1 + q^{-h} \otimes e \\ 
\Delta (f) &=& f \otimes q^{h} + 1 \otimes f \\ 
\end{array} \]

Let us look at the action of $E_j^R $ on a standard monomial: 
\[ E_j^R (z_{1,1}^{d_{1,1}} z_{2,1}^{d_{2,1}} \ldots z_{1,j}^{d_{1,j}} \ldots 
z_{m,j}^{d_{m,j}} z_{1,j+1}^{d_{1,j+1}} \ldots z_{m,j+1}^{d_{m,j+1}} 
\ldots z_{m-1,n}^{d_{m-1,n}} z_{m,n}^{d_{m,n}} ) \]
It is clear the general non-standard term generated will be:
\[ z_{1,1}^{d_{1,1}} z_{2,1}^{d_{2,1}} \ldots z_{1,j}^{d_{1,j}} \ldots 
z_{m,j}^{d_{m,j}} z_{1,j+1}^{d_{1,j+1}} \ldots z_{i-1,j+1}^{d_{i-1,j+1}} 
z_{i,j+1}^a  \cdot z_{i,j} \cdot z_{i,j+1}^b z_{i+1,j+1}^{d_{i+1,j+1}} \ldots  
z_{m,j+1}^{d_{m,j+1}} \ldots z_{m-1,n}^{d_{m-1,n}} z_{m,n}^{d_{m,n}}  \]
where $a+1+b=d_{i,j+1}$. This term straightens to:
\[ q^{a+d_{i+1,j}+\ldots d_{m,j}} 
z_{1,1}^{d_{1,1}} z_{2,1}^{d_{2,1}} \ldots z_{1,j}^{d_{1,j}} \ldots 
z_{i,j}^{d_{i,j}+1} \ldots  z_{m,j}^{d_{m,j}} z_{1,j+1}^{d_{1,j+1}} \ldots 
z_{i-1,j+1}^{d_{i-1,j+1}} z_{i,j+1}^{d_{i,j+1}-1} z_{i+1,j+1}^{d_{i+1,j+1}} 
\ldots  z_{m,j+1}^{d_{m,j+1}} \ldots z_{m-1,n}^{d_{m-1,n}} z_{m,n}^{d_{m,n}}  \]
The Hopf constant will be 
\[ q^{d_{1,j+1}+\ldots +d_{i-1,j+1}+a-d_{1,j}-\ldots -d_{m,j}} \]
Thus the total index $n_{i,j}^R $ is:
\[ \begin{array}{rcl}
n_{i,j,a}^R &=&  a+d_{i+1,j}+ \ldots +d_{m,j}+d_{1,j+1}+\ldots +d_{i-1,j+1}+a-d_{1,j}-\ldots -d_{m,j} \\
&=& (d_{1,j+1}-d_{1,j})+\ldots + (d_{i-1,j+1}-d_{i-1,j})+(2a-d_{i,j}) \\ \end{array} \]
We may abbreviate all this by assuming that $d$ is an $m\times n$-matrix of 
non-negative degrees $d_{i,j}$ and $z^d =\prod_{i,j} z_{i,j}^{d_{i,j}}$. 
Let $\kappa_{ij}^R $ stand for the matrix of all zeros except in the $(i,j)$
and the $(i,j+1)$ positions, where it is $1$ and $-1$ respectively. We may 
thus write:
\[ E_j^R (z^d )=\sum_{i: d+\kappa_{ij}^R \geq 0} \sum_{a=0}^{d_{i,j+1}-1} 
q^{n_{ija}^R} z^{d+\kappa_{ij}^R } \]

Let us look at the easier left action:
\[ \begin{array}{rcl}
  E_i^L (z_{i',j'})&=& \left\{ \begin{array}{lr} 
            z_{i,j'} & \: \: \mbox{ if } i'=i+1 \\
            0   & \mbox{otherwise} \end{array} \right. \\ 
&& \\
  F_i^R (z_{i',j'})&=& \left\{ \begin{array}{lr} 
            z_{i+1,j'} & \: \: \mbox{ if } i'=i \\
            0   & \mbox{otherwise} \end{array} \right. \\ 
\end{array} \]

Let us define $\kappa_{ij}^L $ as the zero matrix, except in the 
$(i,j)$ and $(i+1,j)$-th positions, where it is $1$ and $-1$ respectively, and
\[ \begin{array}{rcl}
n_{i,j,a}^L &=& (d_{i+1,1}-d_{i,1})+\ldots + (d_{i+1,j-1}-d_{i,j-1})+(2a-d_{i,j}) \\ \end{array} \]
In this notation, we see that
\[ E_i^L (z^d )=\sum_{j: d+\kappa_{ij}^L \geq 0}\sum_{a=0}^{d_{i+1,j}-1} 
 q^{n_{ija}^L} z^{d+\kappa_{ij}^L } \]

We are thus led to the following two observations:
\begin{itemize}
\item[(i)] The formulae for the right action are obtained from the left by 
transposing the the first and second indices in $z,d$ etc., i.e., the roles 
of the left and the right action are completely analogous.
\item[(ii)] The left action, i.e., the action of $U_q (gl_m )$ on $S_k (X)$ 
matches the standard action of $U_q (gl_m )$ on the module $\oplus_a 
\otimes_j Sym^{a_j } (\C^m )$ where $\sum_j a_j =k$. Indeed, given a monomial
$z^d $, we read it column-wise to get each component of the $n$-way tensor
product.
\end{itemize}

Since the crystal base for $\otimes Sym^a (\C^m )$ is well understood, we are 
led to the following proposition:

\begin{prop}
A crystal basis for $S^k (X)$ is the collection of monomials $\{ z^d | 
d\geq 0, \sum_{i,j} d_{ij} =k \} $. 
\end{prop}

For a $m\times n$-matrix $d$ of non-negative integers, define the (sym.) 
{\bf left word}
$SLW(d)$ as the $i$-indices of all elements $(i,k)\in d$, {\bf repeated $d_{ik}$ times}, read 
top to bottom within a column, reading the columns left to right. 
Similarly, define the {\bf right 
word} $SRW(d)$ as the $k$-indices of all elements $(i,k) \in d$ {\bf repeated
$d_{ik}$ times}, 
read left to right within a row, reading the rows from top to bottom. 
For a word $w$, let $rs(w)$ be 
the Robinson-Schenstead tableau associated with $w$, when read from
left to right. Define the {\bf left tableau} $SLT(d)=rs(SLW(d))$ and 
the {\bf right tableau} as $SRT(d)=rs(SRW(d))$.

\begin{ex}
Let $m=3$ and $n=4$ and let $d$ be as given below: 
\[ d=\begin{array}{c} 
\young(1002,0201,3110) \: \, \hspace*{1cm}  SLW(d)=13332233112 \: \,
\hspace*{1cm}   SRW(d)= 14422411123 \\
SLT(d)=\young(111233,223,33) \hspace*{1cm} SRT(d)=\young(111123,224,44) \\
\end{array}
\]
\end{ex}

\noindent
We now obtain the result in \cite{danil} for the $Sym$-case.

\begin{prop}
For any $z^d $, element of the crystal basis for 
$S^k (\C^{mn})$ as above, we have:
\begin{itemize}
\item If $\stackrel{\sim}{E_i^L}(z^d )= z^e $ then
$\stackrel{\sim}{e_i^T} (SLT(d))=SLT(e)$. 
\item If $\stackrel{\sim}{E_k^R}(z^d )= z^e $ then
$\stackrel{\sim}{e_k^T} (SRT(d))=SRT(e)$. 
\end{itemize}
A similar assertion holds for the $\stackrel{\sim}{F}$-operators.
\end{prop}

\section{The $3$-column conditions and its verification for the 
$U_q (gl_2 ) \otimes U_q (gl_2 )$ case} \label{sec:22}

Motivated by the construction of the $S(X)$, we define the general algebra 
${\cal F} (X)$ as follows. The generators ${\cal C}$ of the algebra are 
crystal bases of 
each $\wedge^k (\C^{mn})$ indexed by (strict) 
column-tableaus with entries in $[mn]$. Let $T({\cal C})$ be the tensor-algbera
with generators ${\cal C}$. It is clear that $T({\cal C}$ is a $U_q (gl_m ) 
\otimes U_q (gl_n )$-module. Next, given columns $c,c'$, we recall the 
relationship $c\leq c'$ to mean that (i) $|c|\geq |c'|$, and (ii) $c(i)\leq 
c'(i)$ for $i=1,\ldots , |c'|$. In other words, $c'$ may follow $c$ in a 
semi-standard tableau. We use a particular family $\psi_{a,b}$ and 
define the straightening laws:
\[ c\cdot c' =\left\{ \begin{array}{ll}
 \sum_i \alpha_i \cdot c_i \cdot c'_i & \mbox{\hspace*{0.2cm} 
if $|c|\geq |c'|$ and 
$c \not \leq c'$}  \\
 0   & \mbox{\hspace*{0.2cm} if $|c|<|c'|$} 
\end{array} \right. \]
where $c_i \leq c'_i $ for all $i$.
We call ${\cal S}$ as the double-sided ideal generated by the above as the 
straightening relations. Note that ${\cal S}\subseteq T({\cal C})$ is a 
$U_q (gl_m ) \otimes U_q (gl_n )$-module. We define ${\cal F}(X)$ as the 
quotient $T({\cal C})/{\cal S}$. Recall that $SS(\lambda ,[mn])$ is the 
collection of all semi-standard tableau of shape $\lambda $ with entries in 
$[mn]$. Note that ${\cal F}(X)$ has a natural grading and one may hope that:
\[ {\cal F} (X)^d \stackrel{?}{\equiv} \sum_{|\lambda |=d} SS(\lambda ,[mn]) \]

By Bergman's {\em diamond lemma}, the hope above boils down to verifying the 
following $3$-column case. Let $Im_{r,s} \subseteq 
\wedge^r \otimes \wedge^s $ be the image of $\psi_{r,s}$. For $a\geq b \geq 
c>0$ let 
\[ \begin{array}{rcl}
{\cal S}_{a,b,*}&=& Im_{a,b} \otimes \wedge^c  \subseteq \wedge^a \otimes \wedge^b \otimes \wedge^c \\
{\cal S}_{*,b,c}&=& \wedge^a \otimes Im_{b,c} \subseteq \wedge^a \otimes \wedge^b \otimes \wedge^c \\
{\cal S}_{a,b,c}&=& {\cal S}_{*,b,c}+ {\cal S}_{a,b,*} \\
\end{array} \]
Note that ${\cal S}_{a,b,c}$ is a sub-module of $\wedge^a \otimes \wedge^b 
\otimes \wedge^c $. 

For $a,b,c$ as above, let $[a,b,c]'$ denote the $3$-column shape with column
lengths $a,b,c$. We have the obvious:
\begin{prop} \label{prop:3col}
If $dim(\wedge^a \otimes \wedge^b \otimes \wedge^c /{\cal S}_{a,b,c}
)=dim(V_{([a,b,c]')} (\C^{mn}))$ for all $a,b,c$ then 
\[ {\cal F} (X)^d \equiv \sum_{|\lambda |=d} SS(\lambda ,[mn]) \]
\end{prop}

We see that Proposition \ref{prop:3col} brings the problem down to verifying
certain properties of $\psi $ on a finite set of modules and which may even 
be done on a computer. We have done 
precisely this for the $m=n=2$ case. {\em However, we construct special 
$\psi $'s which are not the same as those in Section \ref{sec:psi}}. 

Define $\psi_2 ,\psi'_2 ,\psi_3 $ and $\psi'_3 $ as follows:
\[ \begin{array}{rcl} \vspace*{0.2cm} 
\psi_2 (\hspace*{0.1cm} \young(1,2)\hspace*{0.1cm} ) &=& -q \cdot \young(2) \otimes \young(1) +q^2 \cdot \young(1) \otimes \young(2) \\ \vspace*{0.2cm} 
\psi_2 (\hspace*{0.1cm} \young(1,3)\hspace*{0.1cm} ) &=& -q \cdot \young(3) \otimes \young(1) +q^2 \cdot \young(1) \otimes \young(3) \\ \vspace*{0.2cm} 
\psi'_2 (\hspace*{0.1cm} \young(1,2)\hspace*{0.1cm} ) &=& \young(2) \otimes \young(1) -q \cdot \young(1) \otimes \young(2) \\ \vspace*{0.2cm} 
\psi'_2 (\hspace*{0.1cm} \young(1,3)\hspace*{0.1cm} ) &=& \young(3) \otimes \young(1) -q \cdot \young(1) \otimes \young(3) \\ \vspace*{0.2cm} 
\psi_3 (\hspace*{0.1cm} \young(1,2,3)\hspace*{0.1cm} ) &=& (q+1)(q^2 +q-1)/(1+q^2 ) \cdot \young(2,3) \otimes \young(1) -(q^3 -2q+1)/(1+q^2 ) 
\cdot \young(1,4) \otimes \young(1)  \\ \vspace*{0.2cm} 
 && -(q^2 +q-1)q \cdot \young(1,3) \otimes \young(2) +(q^2 +q-1)q \cdot \young(1,2) \otimes \young(3) \\ \vspace*{0.2cm} 
\psi'_3 (\hspace*{0.1cm} \young(1,2,3)\hspace*{0.1cm} ) &=& -q^3 (q+1)/(1+q^2) 
\cdot \young(1) \otimes \young(2,3) + q^3 (q-1)/(1+q^2 ) \cdot 
\young(1) \otimes \young(1,4) \\ \vspace*{0.2cm} 
&& + q \cdot \young(2) \otimes \young(1,3) -q \cdot \young(3) \otimes \young(1,2) \\ \vspace*{0.2cm} 
\end{array} \]

Note that there are two highest weight vectors in $\wedge^2 (\C^{2 \times 2})$, 
viz., for the shapes $\yng(2)$ and $\yng(1,1)$, while 
$\wedge^3 (\C^{2 \times 2})$ has only one highest weight vector, viz., 
for the shape $\yng(2,1)$ \hspace*{0.05cm}. One may check that the images 
specified are indeed highest weight vectors. Also check that at $q=1$, the 
maps reduce to the classical ones. The surprising terms are, of course, 
$\young(1) \otimes \young(1,4)$ and its counter-part; both vanish at $q=1$.
$\psi_3 $ and $\psi'_3 $ were chosen so that the following diagram commutes:

\begin{center}
\begin{tabular}{c}
\begindc{\commdiag}[20]
\obj(0,6){$\wedge^{3} $}
\obj(7,6){$\wedge^{2} \otimes \wedge^{1}$}
\obj(0,0){$\wedge^{1} \otimes \wedge^{2}$}
\obj(7,0){$\wedge^{1} \otimes \wedge^{1} \otimes \wedge^1 $}
\mor(1,6)(6,6){$\psi_3  $}
\mor(0,6)(0,0){$\psi'_3  $}
\mor(1,0)(5,0){$id \otimes \psi_2  $}
\mor(7,6)(7,0){$\psi'_2 \otimes id $}
\enddc
\end{tabular}
\end{center}

Note that for the purpose of verifying the $3$-column conditions, each of the 
$\psi $'s may be individually scaled. Of course, $\psi_2 $ on $\young(1,2)$ 
and $\young(1,3)$, may be individually scaled while maintaining $U_q (gl_m )
\otimes U_q (gl_n )$-equivariance, but at the risk of changing the 
$3$-column conditions. The maps $\psi_4 $ and $\psi'_4 $ have no real choice; 
there is only one $U_q (gl_m ) \otimes U_q (gl_n )$-invariant in 
$\wedge^3 \otimes \wedge^1 $ or in $\wedge^1 \otimes \wedge^3 $.

We use the above basic maps to construct $\psi_{a,b}$ for all $a\geq b$, 
viz., $\psi_{1,1},\psi_{2,1}, \psi_{2,2},\psi_{3,1}, \psi_{3,2}$ and 
$\psi_{3,3}$. The actual verification of the $3$-column conditions was done
on a computer. The 10 tuples for $a,b,c$ are 
$111,211,221,222,311,321,322,331,332,333$. 

The condition was checked only for the weight space for the left and right 
weights closest to zero. Thus if the number of boxes were even then the left, 
and right weights were chosen $0,0$, else $1,1$, respectively. Clearly, 
a violation of the condition would mean an additional $U_q (gl_2 ) \otimes 
U_q (gl_2 )$-module than in the classical case and would have a witness at
these weights.

\section{Notes} \label{sec:notes}

The immediate objective is to construct $W_{\lambda } (\C^{mn})$ for all 
$\lambda $. One route is through the $3$-column conditions, which would 
mean the construction of a $U_q (gl_m )\otimes U_q (gl_n )$-equivariant
resolution of $V_{[a,b,c]'}$ perhaps mimicking the {\em Giambelli}-type
resolution of 
Akin \cite{akin} in the classical case. This itself is based on the 
Bernstein-Gelfand-Gelfand resolution \cite{bgg}, a $q$-version of 
which is also available. The trouble, of course, is to construct one which is 
$U_q (gl_m ) \otimes U_q (gl_n )$-equivariant. 

The choice of $\psi $ is of course, critical. Shamefully, even for the $m=n=2$
case, other than $a=1,2 $ and $b=c=1$, there seems to be no explanation
of why the $3$-column condition holds. How does the ``commutation'' condition
on $\psi ,\psi'$ actually translate into proofs of the $3$-column conditions?

\end{document}